\documentclass[11pt]{article}
\usepackage{amsmath,amssymb,amsfonts,amsthm,amscd,graphicx,psfrag,epsfig}
\usepackage{color}
\definecolor{blue}{rgb}{0,0,0.7}
\definecolor{red}{rgb}{0.75, 0, 0}
\usepackage{stmaryrd}
\usepackage[titletoc,toc]{appendix}
\usepackage{hyperref}
\hypersetup{colorlinks, linkcolor=blue, citecolor=cyan}
\usepackage[all]{xy}           
\usepackage{marginnote}  
\usepackage{float}
\usepackage{soul}

\newtheorem{theorem}{Theorem}[section]
\newtheorem{lemma}[theorem]{Lemma}
\newtheorem{proposition}[theorem]{Proposition}
\newtheorem{corollary}[theorem]{Corollary}
\newtheorem{conjecture}[theorem]{Conjecture}
\newtheorem{definition}[theorem]{Definition}

\newcommand{\bpf}{\begin{proof}}
\newcommand{\epf}{\end{proof}}
\newcommand{\bs}{\begin{split}}
\newcommand{\es}{\end{split}}
\newcommand{\be}{\begin{equation}}
\newcommand{\ee}{\end{equation}}
\newcommand{\bt}{\begin{theorem}}
\newcommand{\et}{\end{theorem}}
\newcommand{\bd}{\begin{definition}}
\newcommand{\ed}{\end{definition}}
\newcommand{\bp}{\begin{proposition}}
\newcommand{\ep}{\end{proposition}}
\newcommand{\bl}{\begin{lemma}}
\newcommand{\el}{\end{lemma}}
\newcommand{\bc}{\begin{corollary}}
\newcommand{\ec}{\end{corollary}}
\newcommand{\bcon}{\begin{conjecture}}
\newcommand{\econ}{\end{conjecture}}
\newcommand{\la}{\label}
\newcommand{\Z}{{\mathbb Z}}
\newcommand{\R}{{\mathbb R}}
\newcommand{\Q}{{\mathbb Q}}
\newcommand{\C}{{\mathbb C}}
\newcommand{\G}{{\rm G}}
\newcommand{\B}{{\rm B}}
\newcommand{\Id}{{\rm Id}}
\newcommand{\End}{{\rm End}}
\newcommand{\Hom}{{\rm Hom}}
\newcommand{\hra}{\hookrightarrow}
\newcommand{\lra}{\longrightarrow}
\newcommand{\lms}{\longmapsto}

\newcommand{\gr}{{\rm gr}}

\begin{document}

\begin{titlepage}
\title{The Galois group of the category of mixed Hodge -Tate structures}

\author{Alexander Goncharov, Guangyu Zhu}
\date{\it To Alexander Beilinson, for his 60th birthday }
 \end{titlepage}
\maketitle
\tableofcontents

\begin{abstract}

The category ${\rm MHT}_\Q$ of mixed Hodge-Tate structures over $\Q$ is a mixed Tate category. Thanks to the Tannakian formalism  it is equivalent to the category
of graded comodules over a commutative graded  Hopf algebra ${\cal H}_\bullet = \oplus_{n=0}^\infty {\cal H}_n$ over $\Q$.

Since the category  ${\rm MHT}_\Q$ has  homological dimension one,   ${\cal H}_\bullet$ is isomorphic to the commutative graded Hopf algebra provided by the
tensor algebra of the graded vector space given by the  sum  of ${\rm Ext}_{{\rm MHT}_\Q}^1(\Q(0), \Q(n)) = \C/(2\pi i)^n\Q$ over $n>0$.
However this isomorphism is not natural in any sense, e.g. does not exist in families.

We give a  natural  construction of the Hopf algebra ${\cal H}_\bullet$. Namely,   let   $\C^*_\Q:=\C^* \otimes \Q$. Set 
\[{\cal A}_\bullet(\C):=   \Q \oplus \bigoplus_{n=1}^\infty\C_\Q^* \otimes_\Q \C^{\otimes n-1}.
\]
We provide it with a commutative graded Hopf algebra structure, such that  ${\cal H}_\bullet  = {\cal A}_\bullet(\C)$. This implies 
  another construction of the big period map ${\cal H}_n \lra \C_\Q^* \otimes \C$ from \cite{G96}, \cite{G15}.

Generalizing this, we introduce a notion of a Tate dg-algebra $(R, k(1))$, and assign to it  a Hopf dg-algebra ${\cal A}_\bullet(R)$.   

For example, the Tate algebra $(\C, 2\pi i \Q)$  gives rise to the Hopf algebra ${\cal A}_\bullet(\C)$. 

Another  example of a Tate dg-algebra  $(\Omega_X^\bullet, 2\pi i\Q)$ is provided by the holomorphic de Rham complex $\Omega_X^\bullet$ of a complex manifold $X$.  
The sheaf of Hopf dg-algebras ${\cal A}_\bullet(\Omega_X^\bullet)$    
describes a dg-model of the derived category of   variations of Hodge-Tate structures on $X$. The cobar complex  of 
${\cal A}_\bullet(\Omega_X^\bullet)$ is  a dg-model 
for the rational Deligne cohomology of $X$.  

We consider a variant of our construction which starting from  Fontaine's period rings  $\B_{\rm crys}$ /  $\B_{\rm st}$ produces  
graded / dg Hopf   algebras which we 
  relate   to the p-adic Hodge theory. 
\end{abstract}

\section{Introduction and main constructions}

\label{intro}

\noindent



\subsection{Summary} \la{SECC1.1a}

\paragraph{Hopf dg-algebra ${\cal A}_\bullet(R)$ arising from a Tate dg-algebra $R$.} 

  Let $R$ be  a differential graded (dg) algebra, possibly non-commutative,  over a field $k$, with a differential $d$.

A {\it Tate line} $k(1)$ in $R$ is a 1-dimensional $k$-vector subspace of  the center of $R$, such that:

 \begin{itemize}

 \item The  non zero elements in $k(1)$ are invertible, and lie in  the degree $0$  part  of $R$.

 \item
The line $k(1)$   consists of the $d$-constants: $dk(1)=0$.

 \end{itemize}

  \bd A Tate dg-algebra is    a dg-algebra $R$ equipped with   a  Tate line $k(1)$.
 \ed

 Here are important examples of  Tate algebras / sheaves of algebras. In both examples $k=\Q$.

 \begin{enumerate}

 \item $R=\C$, the field of complex numbers, considered as a $\Q$-algebra, with $\Q(1) :=2\pi i\Q$.

  \item $R= \Omega^\bullet_X$, the de Rham complex of sheaves on a  complex manifold $X$, with  $\Q(1) :=2\pi i\Q$.

  Variant: the logarithmic de Rham complex $\Omega^\bullet_{\rm log}$ on a complex  variety,   $\Q(1) :=2\pi i\Q$.

   \end{enumerate}

 Given   a {Tate dg-algebra} $R$, we introduce     a Hopf dg-algebra $\mathcal{A}_\bullet(R)$ with a differential ${\cal D}$.
  It is always non-cocommutative, and commutative if and only if the algebra $R$ is commutative.
  
  The  $k$-vector space $\mathcal{A}_\bullet(R)$  is  defined as follows: 
\be
\begin{split}
&{\cal A}_\bullet(R)=  \bigoplus_{n=0}^\infty{\cal A}_n(R), ~~~~
 {\cal A}_0:=k, ~~ {\cal A}_n:=   \overline  R \otimes \underbrace{{R} \otimes
\ldots \otimes {R}}_{\text{ $n-1$ factors }}.\\
\end{split}
\ee  
It has a  grading $|\cdot |$ induced by the grading of $R$. The differential  ${\cal D}$  is defined  inductively:
\be \la{DIFFDDazz}
{\cal D}(\overline x_1 \otimes x_2 \otimes \ldots \otimes x_n):= d\overline x_1 \cdot  x_2 \otimes (x_3 \otimes \ldots \otimes x_n) + (-1)^{|x_1|} \overline x_1 \otimes {\cal D}(x_2 \otimes \ldots \otimes x_n).
\ee
 
 The coproduct is induced by the deconcatenation map. 
 
 The definition of the product is rather non-trivial, see Section \ref{SEC2.1as}. 

 The Hopf dg-algebra $\mathcal{A}_\bullet(R)$ provides   the  tensor dg-category   
 of dg-comodules over the Hopf dg-algebra $\mathcal{A}_\bullet(R)$.
   Notice that $H^0_{\cal D}(\mathcal{A}_\bullet(R))$ is a graded Hopf algebra, whose grading is induced by the one on $\mathcal{A}_\bullet(R)$.
It provides the abelian tensor category  
 of graded $H^0_{\cal D}(\mathcal{A}_\bullet(R))$-comodules.
If  the grading on $R$ is trivial,  then $H^0_{\cal D}(\mathcal{A}_\bullet(R)) = \mathcal{A}_\bullet(R)$.

 \paragraph{Applications    to Hodge theory.} These constructions have several   applications.

  \begin{enumerate}
  \item   {\it Mixed rational Hodge-Tate structures.} Let $R=\C$,
  $\Q(1) = 2\pi i \Q$. Then the   Hopf algebra $\mathcal{A}_\bullet(\C)$ is \underline{canonically} isomorphic to the
  Tannakian Galois Hopf algebra of the category  of mixed  Hodge-Tate structures over $\Q$. Therefore the category
 of graded $\mathcal{A}_\bullet(\C)$-comodules is canonically equivalent  to the category of mixed $\Q$-Hodge-Tate structures.

 \item    {\it Complexes of variations of rational Hodge-Tate structures.} Let $R=\Omega^\bullet_X$ be the sheaf of commutative dg-algebras given by the
 de Rham complex of a complex manifold $X$, with
  $\Q(1) = 2\pi i \Q \subset {\cal O}_X$. Our construction provides  a   Hopf dg-algebra  $\mathcal{A}_\bullet(\Omega^\bullet_X)$.

  The category   of dg-comodules over   $\mathcal{A}_\bullet(\Omega^\bullet_X)$  is a dg-model of the derived category of complexes of variations of mixed   $\Q$-Hodge-Tate structures  on $X$.

The weight $n$ part   of the cobar complex ${\rm Cobar}_\bullet(\mathcal{A}_\bullet(\Omega^\bullet_X))$ of the Hopf dg-algebra $\mathcal{A}_\bullet(\Omega^\bullet_X)$  is quasi-isomorphic to the  weight $n$ rational Deligne complex on $X$.
  So the commutative dg-algebra ${\rm Cobar}_\bullet(\mathcal{A}_\bullet(\Omega^\bullet_X))$
  is a commutative dg-model for the Deligne cohomology of $X$.

 \item   {\it Variations of rational Hodge-Tate structures.} The graded   Hopf algebra     $H^0_{\cal D}(\mathcal{A}_\bullet(\Omega^\bullet_X))$   is   the Tannakian   Galois Hopf algebra\footnote{An alternative name for this Hopf algebra is the {\it fundamental Hopf algebra}. It plays a role of the fundamental group of a space, and we think about the category as a category of local systems on this space.}   of the category of variations of $\Q$-Hodge-Tate structures on $X$.   So the latter is canonically equivalent to the category of graded
  $H^0_{\cal D}(\mathcal{A}_\bullet(\Omega^\bullet_X))$-comodules.  Therefore the commutative dg-algebra ${\rm Cobar}_\bullet(H_{\cal D}^0(\mathcal{A}_\bullet(\Omega^\bullet_X)))$
  is another  commutative dg-model for the Deligne cohomology of $X$.

 \item   {\it Log-analogs.} If $X$ is a smooth complex variety, there are similar stories   to 2) and 3), obtained by replacing the de Rham complex by the logarithmic de Rham complex $\Omega_{\rm log}^\bullet$.
  \end{enumerate}

\paragraph{A  variant: Tate $\varphi$-algebras and related  Hopf algebras.} Let $k$ be a field. A {\it Tate $\varphi$-algebra} $R$ is    a graded $k$-algebra  $R$   with an invertible
 Tate line  $k(1)  \subset R_1$:
$$
R = \bigoplus_{n \in \Z} R_n, ~~~~k(1)  \subset R_1. 
$$
 Let $k(n) := k(1)^{\otimes n}$, and  
 set
$
\overline {R_n}:= R_n/k(n). 
$
Consider a $k$-vector space
\be
\begin{split}
&{\cal A}^\varphi_\bullet(R)=  \bigoplus_{n=0}^\infty{\cal A}^\varphi_n(R), ~~~~
 {\cal A}^\varphi_0:=k, ~~ {\cal A}^\varphi_n:=   \overline  {R_1} \otimes \underbrace{{{R_1}} \otimes
\ldots \otimes {{R_1}}}_{\text{ $n-1$ factors }}.\\
\end{split}
\ee
We equip it  with a  graded Hopf algebra structure, see Section \ref{sec5.1}.

\paragraph{  An application to $p$-adic crystalline Galois representations.}  Let $\G_{\Q_p}:= {\rm Gal}(\overline \Q_p/\Q_p)$. Recall   Fontaine's algebra  ${\rm B}^+_{\rm dR}$ over $\Q_p$.  It is 
 equipped with a  $\G_{\Q_p}$-action. It has    a canonical
 $\G_{\Q_p}$-submodule of rank $1$ over  $\Q_p$, generated  by Fontaine's $t$, the p-adic analog of $2 \pi i$:
 $$
\Q_p(1) = \Q_p t \subset {\rm B}^+_{\rm dR}.
$$

Then ${\rm B}_{\rm dR}:= {\rm B}^+_{\rm dR}[t^{-1}]$ is a field. 
There is a $\G_{\Q_p}$-invariant subalgebra ${\rm B}_{\rm crys} \subset {\rm B}_{\rm dR}$. 
It is equipped with an action of the Frobenius $\varphi$. 
Denote by ${\rm B}_{\rm crys}^{\varphi =p^n}$ the subspace where $\varphi$ acts with the eigenvalue $p^n$. Then there is the fundamental exact sequence
$$
0 \lra \Q_p \lra {\rm B}_{\rm crys}^{\varphi =1} \lra {\rm B}_{\rm dR}/{\rm B}^+_{\rm dR} \lra 0.
$$
Multiplying by $t^n$, it implies that for each $n \in \Z$ we have an exact sequence:
\be \la{FES}
0 \lra \Q_pt^n \lra {\rm B}_{\rm crys}^{\varphi =p^n} \lra {\rm B}_{\rm dR}/t^n{\rm B}^+_{\rm dR} \lra 0.
\ee
Let us set 
$$
R_n:= {\rm B}_{\rm crys}^{\varphi =p^n}.
$$
Then there is a subring   ${\rm R}_{\rm crys} \subset \B_{\rm crys}$ with a structure of the Tate $\varphi$-algebra over $\Q_p$:
\be \la{TAQP}
  {\rm R}_{\rm crys} := \bigoplus_{n \in \Z} {\rm B}_{\rm crys}^{\varphi =p^n}, ~~~~ \Q_p(1)  \subset {\rm B}_{\rm crys}^{\varphi =p}.
\ee
Thanks to the  exact sequence (\ref{FES})   we have
$
\overline {R_n} = {\rm B}_{\rm dR}/t^n{\rm B}^+_{\rm dR}.
$

The  Tate $\varphi$-algebra  $({\rm R}_{\rm crys}, \Q_p(1))$ gives rise to a graded Hopf algebra ${\cal A}^\varphi_\bullet( {\rm R}_{\rm crys})$:
\be \la{ACRYS}
\begin{split}
&{\cal A}^\varphi_\bullet({\rm R}_{\rm crys}):= \Q_p\oplus \bigoplus_{n=1}^\infty     
  {\rm B}_{\rm dR}/t{\rm B}^+_{\rm dR} \otimes_{\Q_p} \underbrace{{\rm B}_{\rm crys}^{\varphi =p}
 \otimes _{\Q_p}   \ldots \otimes _{\Q_p}  {\rm B}_{\rm crys}^{\varphi =p}}_{\text{ $n-1$ factors }}.\\
\end{split}
\ee
We define in  Section \ref{sec5.3} a functor from the   category of crystalline mixed Tate p-adic $\G_{\Q_p}$-representations  to the category 
of graded comodules over the Hopf algebra ${\cal A}^\varphi_\bullet({{\rm R}_{\rm crys}})$. 

One of its applications    is an explicit  description of 
the p-adic regulator map on the category of mixed Tate motives. 

In particular, let $F$ be a number field, and $F_v$ the completion of $F$ at a  
non-archimidean valuation $v$ on $F$. Then for  any integer $n>1$ we get a p-adic regulator map
\be \la{Bh}
\begin{split}
&R_{F_v}: K_{2n-1}(F)\lra  {\rm B}_{\rm dR}/t^n{\rm B}^+_{\rm dR}.\\
\end{split}
\ee  
It is the non-archimedean analog of the regulator map assigned to a given embedding $\sigma: F \hra \C$:
\be \la{Bh}
\begin{split}
&R_{\C, \sigma}: K_{2n-1}(F)\lra  \C/(2\pi i \Q)^n\C.\\
\end{split}
\ee

A   general construction of the p-adic  regulator maps, with the target given by the   p-adic syntomic cohomology,
 is well known, see \cite{NN}, \cite{Ni1}, \cite{Ni2} and references there. 
Our construction   allows to construct p-adic regulator maps explicitly using p-adic polylogarithms and their generalizations, see \cite{GZ},
 just  like the one in the complex case in \cite{G15}.
 
 \paragraph{A  Hopf dg-algebra arising from ${\rm B}_{\rm st}$ and semi-stable mixed Tate Galois modules.} We develop a similar story starting with
   Fontaine's period ring $\B_{\rm st}$. First, we use the monodromy  operator $N$ to produce a dg-algebra   ${\Bbb B}^\bullet_{\rm st} := \B_{\rm st} \stackrel{N}{\lra} \B_{\rm st} $   with the differential $N$. 
   We define a dg-subalgebra ${\rm R}^\bullet_{\rm st}\subset {\Bbb B}^\bullet_{\rm st}$. We use the dg-algebra $({\rm R}^\bullet_{\rm st}, N)$ in a way similar to the  use of the de Rham dg-algebra  
   $(\Omega^\bullet(X), d)$ of a regular variety $X$. Namely, it 
  gives rise to a Hopf dg-algebra ${\cal A}^\varphi_\bullet({\rm R}^\bullet_{\rm st})$. The category of dg-comodules over  
 ${\cal A}^\varphi_\bullet({\rm R}^\bullet_{\rm st})$ is a dg-version of the derived category of  mixed Tate semi-stable p-adic Galois representations.  
This will be discussed in details in \cite{GZ}. 

\vskip 3mm

  In Sections \ref{SECC1.2}-\ref{SECC1.4} we discuss in detail
  the crucial case of $\Q$-mixed Hodge-Tate structures.

\subsection{Main example: the category of mixed Hodge-Tate structures}\la{SECC1.2}

Below all tensor products are over the field $\Q$, and $\C^*$ denotes
the $\Q$-vector space $\C^* \otimes \Q$.

    A {\it mixed $\Q$-Hodge-Tate structure}  is a $\Q$-vector space $V$  equipped with a weight filtration $W_{\bullet}$,
 and a Hodge filtration $F^{\bullet}$ of its complexification $V_\C$,
such that $\gr_{m}^WV=0$ for odd  $m$, and  the filtration $F^{\bullet}$
and its conjugate $\overline F^{\bullet}$ induce on $\gr_{2n}^WV$ a pure
weight $2n$  Hodge-Tate structure:
$$
\gr^W_{2n}V_\C =  F_{(2n)}^{n}\cap \overline F_{(2n)}^{n}.
$$
Here $F_{(2n)}^{\bullet}$ is the filtration on $\gr_{2n}^WV_\C$ induced by
$F^{\bullet}$, and similarly $\overline F_{(2n)}^{\bullet}$.
Equivalently, $\gr^W_{2n}V$ is  a direct sum of the Hodge-Tate structures $\Q(-n)$.

The category   ${\rm MHT}_\Q$  of mixed   Hodge-Tate structures over $\Q$
 is the smallest tensor subcategory of the category of all
mixed  Hodge structures over $\Q$ which contains the Hodge-Tate structures
$\Q(0)$ and $\Q(1)$ and closed under the extensions.

\paragraph{The  Tannakian   Galois group of the category ${\rm MHT}_{\Q}$.} The category ${\rm MHT}_{\Q}$  is a
Tannakian category
 with
a  fiber functor $\omega $ to the category  of $\Q$-vector spaces:
$$
\omega : {\rm MHT}_{\Q} \lra {\rm Vect}_\Q, ~~~~H \lms {\rm gr}^WH.
$$
It assigns to a   mixed Hodge-Tate structure  $H$  the
underlying rational vector space of the associate graded ${\rm gr}^WH$.
The {\it  Tannakian Galois group} $\G_{\rm MHT}$  of the category ${\rm MHT}_{\Q}$ is   a  pro-algebraic group over $\Q$, defined as the group of
automorphisms of the fiber functor respecting the tensor structure:
$$
\G_{\rm MHT}:= {\rm Aut}^\otimes\omega .
$$
Thanks to the Tannakian formalism, the functor $\omega $ induces  an equivalence of
 tensor categories
$$
\omega: {\rm MHT}_{\Q} \stackrel{\sim}{\lra} \mbox{finite dimensional}~\G_{\rm MHT}-{\rm modules}.
$$
Denote by ${\rm HT}_\Q$ the  category of pure Hodge-Tate structures. There are two functors
\be
\begin{split}
&\gr^W: {\rm MHT_{\Q}} \lra {\rm HT}_\Q, ~~ V \lms \gr^WV, \\
&i: {\rm HT}_\Q \hra {\rm MHT}_{\Q}, ~~ V \lms V, ~~~~
\gr^W\circ i = {\rm Id}. \\
\end{split}
\ee
The functor $i$ is the canonical embedding.  The category ${\rm HT}_\Q$  is equivalent to the category of
representations of the multiplicative group ${\Bbb G}_m$. Therefore there are two homomorphisms
$$
s: {\Bbb G}_{m}\lra \G_{\rm MHT}, ~~~~p: \G_{\rm MHT}\lra {\Bbb G}_{m}, ~~~~p\circ s={\rm Id}.
$$
They provide  the  group $\G_{\rm MHT}$ with a   structure of a
semidirect product:
$$
0\lra {\rm U}  \lra \G_{\rm MHT} \stackrel{p}{\lra} {\Bbb G}_{m }\lra 0, ~~~~{\rm U} :={\rm Ker}~p.
$$

The action of  $s({\Bbb G}_{m})$ provides   the Lie algebra of the group ${\rm U} $  with
 a structure of a graded pro-Lie algebra over $\Q$. We denote it by  ${\rm L}_{\bullet}$ and call  the Hodge-Tate Lie algebra.

 Equivalently, ${\rm L}_{\bullet}$ can be viewed as   a Lie algebra in the category $ {\rm HT}_\Q$.

 \bd
The Tannakian Galois Hopf algebra   ${\cal H}_\bullet$ of the category ${\rm MHT}_\Q$, also called  the Hodge-Tate  Hopf algebra,  is the commutative Hopf algebra of  regular functions on the pro-algebraic group ${\rm U}$, with the grading provided by the action of $s({\Bbb G}_m)$:
 $$
 {\cal H}_\bullet = {\cal O}({\rm U}).
  $$
  \ed
The   Hopf algebra   ${\cal H}_\bullet$  is the graded dual of the universal enveloping algebra
 of the graded Lie algebra ${\rm L}_{\bullet}$. So it is graded by non-negative integers:
 $$
 {\cal H}_\bullet = \bigoplus_{n=0}^\infty {\cal H}_n.
   $$

   The Tannakian definition of the   Hodge-Tate Hopf algebra  ${\cal H}_\bullet$ is  rather abstract.
   There is a direct description of the Hodge-Tate Hopf algebra  ${\cal H}_\bullet$  via framed objects   of the category ${\rm MHT}_\Q$.
    For    convenience of the reader we recall it now borrowing from  \cite{BGSV} and \cite[Section 4]{G96}.

  \paragraph{The Hodge-Tate Hopf algebra  ${\cal H}_\bullet$ of the category ${\rm MHT}_\Q$ via framed objects.}
      An {\it $n$-framing} on a mixed Hodge-Tate structure  $H$ over $\Q$ is a choice of a
nonzero maps
$$
v_0: \Q(0) \to \gr^W_{0}H, ~~~~f^n: \gr^W_{-2n}H \to \Q(n).
$$
 Denote by $(f^n, H, v_0)$  an $n$-framed Hodge-Tate structure. Consider the equivalence
relation $\sim$ on the set of  $n$-framed mixed Hodge-Tate structures induced by the following condition: any  morphism $H_1 \to H_2$ in the category ${\rm MHT}_\Q$ compatible with frames i.e. making the following diagrams commute

~~~~~~~~~~~~~~~~~~~~~~\xymatrix{
\mathbb{Q}(0) \ar[r] \ar[d]^{\Id} & \gr^W_{0} H_1\ar[d] & & \gr^W_{-2n}H_1 \ar[r] \ar[d] & \mathbb{Q}(n) \ar[d]^{\Id}\\
\mathbb{Q}(0) \ar[r]& \gr^W_{0} H_2 & & \gr^W_{-2n}H_2 \ar[r] & \mathbb{Q}(n),}

\noindent gives rise to an equivalence
 $H_1 \sim H_2$.
Let ${\cal H}_{n}$ be the set of equivalence classes.
Then the set  ${\cal H}_{n}$ has a $\Q$-vector space structure defined as follows:
\be
\begin{split}
&(f^n, H, v_0) + (\widetilde f^n, \widetilde H, \widetilde v_0):=
(f^n+ \widetilde f^n, H\oplus \widetilde H, v_0 + \widetilde v_0); \\
&\lambda (f^n, H, v_0) := ( \lambda f^n, H, v_0). \\
\end{split}
\ee
The tensor product of mixed Hodge-Tate structures induces the commutative multiplication
$$
\mu: {\cal H}_{p} \otimes {\cal H}_{q}\to {\cal H}_{p+ q}.
$$
Next, there is   a coproduct
 \begin{equation} \label{GH}
\Delta=\bigoplus_{p+q=n} \Delta_{p,q}: {\cal H}_{n}\to \bigoplus_{p+q = n}
{\cal H}_{p}\otimes
{\cal H}_{q}.
 \end{equation}
Namely, let $(f^n, H, v_0)\in {\cal H}_{n}$.
Choose a basis $\{v^\alpha_p\}$ in ${\rm Hom}(\Q(p), \gr_{-2p}^W H)$
and the dual basis $\{f_\alpha^p\}$ in
${\rm Hom}(\gr_{-2p}^W H, \Q(p))$.    Then the $(p,q)$ component   of the coproduct is given by:
$$
 (f^n, H, v_0)\lms  \sum_\alpha (f^n, H, v^\alpha_p)(-p) \otimes (f_\alpha^p, H, v_0).
$$

\begin{theorem} \label{Theorem 4.1}
 The graded $\Q$-vector space
${\cal H}_\bullet:= \oplus_{n=0}^\infty {\cal H}_n$
is a graded Hopf algebra with the commutative product
$\mu$ and the coproduct $\Delta$.
The fundamental Hopf algebra of the category ${\rm MHT}_\Q$   is canonically isomorphic to the graded Hopf algebra  ${\cal H}_{\bullet}$.
\end{theorem}

From now on we consider the fundamental Hopf algebra of the category ${\rm MHT}_\Q$ as the Hopf algebra  ${\cal H}_{\bullet}$ of framed objects in the category ${\rm MHT}_\Q$.

\vskip 3mm
The next question is what is the structure of the Hopf algebra ${\cal H}_{\bullet}$ as an abstract Hopf algebra. The answer is well known, as explained in the next paragraph.

\paragraph{The structure of the Hodge-Tate Lie algebra ${\rm L}_{\bullet}$.} A.A. Beilinson proved \cite{Be86} that
\be \la{eq7}
{\rm Ext}^{>1}_{{\rm MHT}_\Q}(\Q(0), \Q(n))=0.
\ee
 Thus the only non-zero Ext-groups  between $\Q(0)$ and $\Q(n)$ in this the category are the  following:
\be \la{1}
\begin{split}
&{\rm Hom}_{{\rm MHT}_\Q}(\Q(0), \Q(0)) = \Q,\\
&{\rm Ext}^1_{{\rm MHT}_\Q}(\Q(0), \Q(n)) = \C^*(n-1) := \frac{\C}{(2\pi i)^n\Q}~~\mbox{for $n>0$}.\\
\end{split}
\ee
Therefore the Hodge-Tate Lie algebra ${\rm L}_{\bullet}$ is a free Lie algebra in the category ${\rm HT}_\Q$ generated by
\be \la{11}
\begin{split}
&\bigoplus_{n >0}{\rm Ext}_{{\rm MHT}_\Q}^1(\Q(0), \Q(n))^{\vee}\otimes \Q(n) = \\
&\bigoplus_{n >0}  \C^*(n-1)^{\vee}\otimes \Q(n).\\
\end{split}
\ee

Recall that the tensor algebra ${\rm T}_\bullet(V)$ of any   vector space $V$ has a natural structure of a graded
commutative Hopf algebra.  The coproduct $\Delta$ is given by the deconcatenation map.
So the restricted coproduct $\Delta':= \Delta - (1\otimes {\rm Id} + {\rm Id}\otimes 1)$ looks as follows:
\be \la{RESTRCOP}
\Delta': v_1 \otimes \ldots \otimes v_n \lms   \sum_{p=1}^{n-1} (v_1 \otimes \ldots \otimes v_p) \bigotimes (v_{p+1} \otimes \ldots \otimes v_n).
\ee
The commutative product $\ast$ is given by the shuffle product:
$$
(v_1 \otimes \ldots \otimes v_p) \ast (v_{p+1} \otimes \ldots \otimes v_{p+q}):= \sum_{\sigma} v_{\sigma(1)} \otimes \ldots \otimes v_{\sigma(p+q)}.
 $$
Here the sum over all snuffles $\sigma$ of the sets $\{1, ..., p\}$ and $\{p+1, ..., p+q\}$.
The grading is given by the  tensors  degree.
If $V_\bullet$ is a  graded vector space, then the Hopf algebra ${\rm T}_\bullet(V_\bullet)$ has a grading given by the sum of the grading on $V_\bullet$.

Consider the graded $\Q$-vector space $E_\bullet$ given by the direct sum of
spaces ${\rm Ext}^1_{{\rm MHT}_\Q}(\Q(0), \Q(n))$ placed in the degree $n$:
$$
E_\bullet:= \oplus_{n=1}^\infty \C^*(n-1)[-n].
$$
The commutative graded Hopf algebras ${\cal H}_\bullet$ and
  ${\rm T}_\bullet(E_\bullet)$ are isomorphic:
\be \la{eq1}
{\cal H}_\bullet \stackrel{\sim}{=} {\rm T}_\bullet(E_\bullet).
\ee
However  \underline{there is no natural   isomorphism} (\ref{eq1}). Here is the simplest example demonstrating this.

\paragraph{Example: the weight 2.} The degree $2$ part of  ${\rm T}_\bullet(E_\bullet)$ is given by
$$
{\rm T}_2(E_\bullet) = \C^*(1) \bigoplus \C^* \otimes \C^*.
$$
Any choice of an isomorphism
$$
{\cal H}_2 \stackrel{\sim}{\lra}   \C^*(1) \bigoplus \C^* \otimes \C^*
$$
would imply that any Hodge-Tate structure $(f^2, H, v_0)$ framed by $\Q(0)$ and $\Q(2)$ has a canonical $\C^*(1)$-valued invariant, a period.
However   such a period does not exist in families. For example, consider the variation of framed Hodge-Tate structures
on $\C{\Bbb P}^1 - \{0,1, \infty\}$ corresponding to the  dilogarithm ${\rm Li}_2(z)$. Then, due to the multivalued nature of the dilogarithm,
 there is no way to assign to it a  $\C^*(1)$-valued invariant   depending continuously on
$z\in \C{\Bbb P}^1 - \{0,1, \infty\}$.

\vskip 3mm
In this paper we suggest a natural
explicit construction of the  Hopf algebra ${\cal H}_\bullet$.

\paragraph{The graded $\Q$-vector space    ${\cal A}_\bullet$.}  Let us consider the following
graded $\Q$-vector space:
\be
\begin{split}
&{\cal A}_\bullet:= \bigoplus_{n=0}^\infty{\cal A}_n,\\
&{\cal A}_0:=\Q, ~~{\cal A}_1 = \C^\ast, ~~{\cal A}_2 = \C^\ast\otimes \C, ~~  \ldots ,~~~ {\cal A}_n:=   \C^*\otimes \underbrace{{\C} \otimes
\ldots \otimes {\C}}_{\text{$n-1$ factors }}.\\
\end{split}
\ee

\paragraph{The coalgebra structure  on ${\cal A}_\bullet$.} We are going to define a coproduct map
$$
\Delta: {\cal A}_\bullet \lra  {\cal A}_\bullet \otimes {\cal A}_\bullet.
$$
The restricted coproduct  $\Delta'$ is decomposed into components
$$
\Delta':= \sum_{p, q\geq 1}\Delta_{p,q}, ~~~~ \Delta_{p,q}: {\cal A}_{p+q} \lra {\cal A}_p \otimes {\cal A}_{q}.
$$
The map $\Delta_{p,q}$  is given by applying the exponential map to
 the $(p+1)$-st factor of the tensor product $\C^*\otimes \C \otimes \ldots \otimes \C  $,
and breaking the tensor product into a product of two factors:
\be \la{eq4}
\begin{split}
&\Delta_{p,q}:\ \C^* \otimes  \C \otimes \ldots \otimes \C  \lra
 \C^* \otimes  \underbrace{{\C} \otimes \ldots \otimes {\C}}_{\text{ $p-1$ factors }}
\bigotimes  \C^* \otimes \underbrace{{\C} \otimes
\ldots \otimes {\C}}_{\text{ $q-1$ factors }}. \\
&\Delta_{p,q}:=\underbrace{{\rm Id} \otimes \ldots \otimes {\rm Id}} _{\text{ $p$ factors }}
 \bigotimes  \exp \otimes  \underbrace{{\rm Id} \otimes \ldots \otimes {\rm Id}} _{\text{ $q-1$ factors }}. \\
\end{split}
\ee
The coproduct  is evidently coassociative.
Our first main result  is the following.
\bt \la{Th1.3} The  graded coalgebra ${\cal A}_\bullet$ has a  commutative product, making it into a graded commutative  Hopf algebra, which is
   \underline{canonically} isomorphic to the  Hodge-Tate Hopf algebra ${\cal H}_\bullet$:
$$
{\cal A}_\bullet \cong  {\cal H}_\bullet.
$$
\et
The   construction of the product  is given in Theorem \ref{Th1.4}.

Below we elaborate   two examples,  defining the commutative Hopf algebra structure on the weight two  and three  quotients ${\cal A}_{\leq 2}$ and ${\cal A}_{\leq 3}$ of
the Hopf algebra ${\cal A}_\bullet$.

\paragraph{Example: the weight 2.} Recall that ${\cal A}_2 = \C^* \otimes \C$. One has an exact sequence
\be \la{eq3}
0 \lra \C^*(1) \lra  \C^* \otimes \C \lra \C^* \otimes \C^* \lra 0.
\ee
Here the first map is given by $a \otimes 2\pi i \lms a \otimes 2\pi i$, and the second is $a \otimes b \lms a \otimes e^b$.

Sequence (\ref{eq3}) is obtained by taking the tensor product of $\C^*$ with the exponential sequence
$$
0 \lra \Q(1) \lra    \C \stackrel{\rm exp}{\lra}   \C^* \lra 0.
$$
The exponential sequence does not have a natural splitting,  meaning that its version over a base
$$
0 \lra \Q_X(1) \lra    {\cal O}_X \stackrel{\rm exp}{\lra}    {\cal O}_X^* \lra 0.$$
 does not have a continuous splitting.
Indeed,   a splitting  requires to  choose a branch of the logarithm  $\log(f)$, $f \in  {\cal O}_X^*$.

Let us now describe explicitly  the weight $2$ quotient Hopf algebra
$$
{\cal A}_{\leq 2} = \Q ~\bigoplus~ \C^* ~\bigoplus ~\C^* \otimes \C.
$$
The only non-trivial component of the restricted coproduct $\Delta'$ is given by
\be
\begin{split}
\Delta^\prime: \ \ & \C^* \otimes \C \lra \C^* \otimes \C^*, \\
&A \otimes b\lms A \otimes e^b.\\
 \end{split}
   \ee
 The commutative product $\ast$ is described by the map
\be
\begin{split}
\ast: \ \ &{\rm Sym}^2 \C^* \lra \C^* \otimes \C, \\
&A \ast B =  A \otimes \log B + B \otimes \log A -  \exp\left(\frac{\log A \log B}{2\pi i}\right)\otimes 2\pi i. \\
\end{split}
\ee
Here $\log A$ denotes a branch of the logarithm. Altering  it by $\log A \lms \log A + 2\pi i n$ we do not change the right hand side.
Similarly for $\log B$.

We stress that in contrast with the shuffle product formula for the product on the Hopf algebra ${\rm T}(E_\bullet)$, already on the quotient  ${\cal A}_{\leq 2}$ the product is  rather non-trivial.

\paragraph{Example: the weight 3 quotient Hopf algebra ${\cal A}_{\leq 3}$.} We have
$$
{\cal A}_{\leq 3} = \Q ~\bigoplus ~\C^* ~\bigoplus ~\C^* \otimes \C~\bigoplus ~\C^* \otimes \C\otimes \C.
$$
The coproduct is given by (\ref{eq4}). The non-trivial components of the product are the following.
\be
\begin{split}
\C^* \ast \left(\C^* \otimes \C \right) \lra &\ \C^* \otimes \C \otimes \C;\\
A_1 \ast (B_1 \otimes b_2) \lms &\  A_1 \otimes \log B_1 \otimes b_2 +  B_1 \otimes \log A_1 \otimes b_2 +  B_1 \otimes b_2 \otimes \log A_1 -\\
 &\ \exp \left( \frac{\log A_1 \log B_1}{2\pi i}\right)  \otimes 2\pi i  \otimes b_2
 -B_1 \otimes  \frac{(\log A_1) b_2}{2\pi i}  \otimes 2\pi i
 .\\
 \end{split}
   \ee
   Altering the branches of  logarithms by $\log A_1 \lms \log A_1 + 2\pi i n$ and $\log B_1 \lms \log B_1 + 2\pi i m$ leaves the formulas intact.

\paragraph{A graded commutative Hopf algebra $\mathcal{A}_\bullet$.} Our next goal is to define a commutative product
on the graded space ${\cal A}_\bullet$ which generalizes the above formulas.

Recall that  the space ${\rm T}_\bullet(\C)$ of $\C$ has a commutative shuffle product.

We define another commutative product $m'$ on the graded space ${\rm T}_\bullet(\C)$,
  using the  induction:
  \be \label{EQP}
  \begin{split}
  m'(x_1\otimes\cdots \otimes x_p, y_1\otimes\cdots \otimes y_q)&:= x_1 \otimes m'(x_2\otimes\cdots \otimes x_p, y_1\otimes\cdots \otimes y_q)\\
  &+ y_1 \otimes m'(x_1\otimes\cdots \otimes x_p, y_2\otimes\cdots \otimes y_q)\\
  &- \frac{x_1 y_1}{2\pi i} \otimes 2\pi i \otimes m'(x_2\otimes\cdots \otimes x_p, y_2\otimes\cdots \otimes y_q).
\end{split}
\ee

However the product $m'$ is \underline{not compatible} with the standard  coproduct  on ${\rm T}_\bullet(\C)$.

Let us consider the canonical projection of graded vector spaces
\be \la{PROJJ}
\begin{split}
 {\rm pr}: \ \ &{\rm T}_\bullet(\C) \lra \mathcal{A}_\bullet,\\
{\rm pr}: \ \ &  \C^{\otimes n}\lra   \C^*  \otimes \C^{\otimes n-1},\\
& x_1\otimes x_2\otimes\cdots \otimes x_n \lms e^{x_1} \otimes x_2\otimes\cdots \otimes x_n.\\
\end{split}
\ee

Here is the precise version of our first main result, Theorem \ref{Th1.3}.

\begin{theorem} \la{Th1.4} The product $m'$ on ${\rm T}_\bullet(\C)$ descends under the map ${\rm pr}$ to a  commutative graded product $m$ on $\mathcal{A}_\bullet$.
It is compatible with the coproduct $\Delta$ on $\mathcal{A}_\bullet$. The triple  $(\mathcal{A}_\bullet, m, \Delta)$ is  a connected commutative graded Hopf algebra.

The   Hopf algebra $(\mathcal{A}_\bullet, m, \Delta)$ is canonically isomorphic to the Hodge-Tate Hopf algebra $\mathcal{H}_\bullet$.
\end{theorem}

Our next goal is  an explicit formula for the  product $m$ on ${\rm T}_\bullet(\C)$.

\subsection{Quasi-shuffle product on ${\rm T}_\bullet(\C)$} \la{SECC1.3}

\paragraph{Quasi-shuffle product $\ast$.} We define a {\it quasi-ordered set} as a set $S$ presented as a disjoint   union of subsets $S_i$, together with an order $<$ of the subsets:
\be \la{qs}
S = S_1 \cup \ldots \cup S_n, ~~~~ S_1 < S_2 < \ldots < S_n.
\ee
We use a notation $S = \{S_1, \ldots , S_n\}$ for a quasi-ordered set. We picture it by a collection of distinct points $\{s_1, ..., s_n\}$ on the real line
with the subset $S_i$ attached to     the point $s_i$.

Let us   define  a   {\it quasi-shuffle} of  two quasi-ordered   sets.

Consider a  collection  of points $\{s_1, ..., s_p, t_1, ..., t_q\}$
on the real line such that

\begin{itemize}

\item The points $\{s_1, ..., s_p\}$ are distinct points, whose   order  on the line coincides with their given order. Same about
the points $\{t_1, ..., t_q\}$.

\item Some of the points $s_i$ may coincide with some  of the points $t_j$. We refer to such a pair of point  as  a {\it colliding pair of points} $s_i=t_j$.
We may have several colliding pairs of points.

The rest of the points are referred to as {\it single points}.

\end{itemize}

We assign to each point of such a collection  a subset as follows:

 1) We assign to a single point $s_i$ (respectively $t_j$) the subset $S_i$ (respectively $T_j$).

 2) We assign to a colliding pair of points $s_i=t_j$ the union of subsets $S_i \cup T_j$.

This way we get a new  quasi-ordered set,   called a {\it quasi-shuffle} of the  original ones.

Given quasi-ordered   sets $\{{S}_1, \ldots , {S}_p\}$ and $\{{T}_1, \ldots , {T}_q\}$,  their
{\it quasi-shuffle product} is defined as a signed formal sum of all possible up to an isomorphism\footnote{quasi shuffles $S$ and $S'$ are isomorphic if there is an isomorphism of ordered sets $S \to S'$ respecting the decompositions (\ref{qs}).} quasi-shuffles. Precisesly, we have:

Consider the free abelian group generated by the isomorphism classes of quasi-ordered sets. 

 \bd Quasi-shuffle product of    two disjoint quasi-ordered sets  $\{{S}_1, \ldots , {S}_p\}$ and $\{{T}_1, \ldots , {T}_q\}$ is a formal signed sum
of all their quasi-shuffles,  considered up to an isomorphism,  taken with the sign $(-1)^c$ where $c$ is the number of pairs of colliding points.
\ed
We use the   notation $
\{{S}_1, \ldots , {S}_p\}\ast \{{T}_1, \ldots , {T}_q\}$ for the quasi-shuffle product.
For example:
\be
\begin{split}
\{{S}_1\}\ast \{{T}_1\} = &\ \{{S}_1, {T}_1\} + \{ {T}_1, {S}_1\} - \{{S}_1 \cup  {T}_1\}.\\
&\\
\{{S}_1, {S}_2\}\ast \{{T}_1\} = &\ \{{S}_1, {S}_2, {T}_1\} + \{{S}_1, {T}_1, {S}_2\} + \{{T}_1, {S}_1, {S}_2\} - \{{S}_1,  {S}_2\cup  {T}_1\}- \{{S}_1\cup {T}_1,  {S}_2\}.\\
&\\
\{{S}_1, {S}_2\}\ast \{{T}_1, {T}_2\} =
& \ \ \{{S}_1, {S}_2, {T}_1, {T}_2\}
+ \{{S}_1,  {T}_1, {S}_2, {T}_2\}
+ \{{S}_1,  {T}_1, {T}_2, {S}_2\}\\
&+ \{{T}_1, {S}_1, {S}_2,  {T}_2\}
+ \{{T}_1, {S}_1,   {T}_2, {S}_2\}
+\{{T}_1, {T}_2, {S}_1, {S}_2\}  \\
& - \{{S}_1,  {S}_2\cup {T}_1, {T}_2\} -
 \{{S}_1\cup {T}_1, {S}_2, {T}_2\} -
  \{{S}_1,   {T}_1, {S}_2\cup {T}_2\}   \\
& +  \{{S}_1\cup{T}_1,  {S}_2\cup {T}_2\}.\\
\end{split}
\ee

Similarly one can define a quasi-shuffle product of any finite number of quasi-ordered sets.
It is evidently commutative and associative.

\paragraph{A geometric origin of the quasi-shuffle product.}  Recall  the standard  simplex
$$
\Delta_m:= \{0 \leq t_1 \leq \ldots \leq t_m \leq 1\}, ~~~~t_i \in \R.
$$
Each quasi-shuffle ${\bf q}$ of the  ordered sets $\{1, ..., m\}$ and $\{1, ..., n\}$ determines     a simplex $\Delta_{\bf q}$ in
the standard simplicial decomposition   of the product $\Delta_m \times \Delta_n$:
$$
\Delta_{\bf q} \subset \Delta_m \times \Delta_n.
$$

Let ${\rm cod}_{\bf q} $ be the codimension of the simplex $\Delta_{\bf q}$ in $ \Delta_m \times \Delta_n$.
Then the product of the two simplices is presented by a quasi-shuffle product  inclusion - exclusion formula:
\be
\label{DEC}
\Delta_m \times \Delta_n = \sum_{\bf q}(-1)^{{\rm cod}_{\bf q}}\Delta_{\bf q}.
\ee
For example, the product of two closed intervals is a union of two triangles minus the diagonal:
\be \la{QS1}
\{0\leq t_1, t_2 \leq 1\} = \{0 \leq t_1 \leq  t_2 \leq 1\} \cup \{0 \leq t_2 \leq  t_1 \leq 1\} - \{0 \leq t_1 = t_2 \leq 1\}.
\ee
We conclude that the quasi-shuffle product reflects the decomposition (\ref{DEC}).

\paragraph{An explicit formula for the product.} We define the product
$$
m'(x_1\otimes\cdots \otimes x_p, y_1\otimes\cdots \otimes y_q)
$$
by taking the quasi-shuffle product of the ordered sets $\{x_1, ..., x_p\}\ast \{y_1, ..., y_q\}$ and  assigning to it an element of $\C^{\otimes p+q}$ as follows.
For a given quasi-shuffle:

1. We assign to each single point $x_i$ (respectively $y_i$) the element $x_i \in \C$ (respectively $y_i$);

2. We assign to each colliding pair of points $x_i=y_j$ the element $ \frac{x_iy_i}{2\pi i}\otimes 2\pi i$.

Then we take the signed sum over all quasi-shuffles.

For example, we have:
\be \la{QSPFa}
\begin{split}
&m'(x_1, x_2) =\  x_1\otimes x_2 + x_2\otimes x_1 -  \frac{x_1x_2}{2\pi i}\otimes 2\pi i.\\
&\\
&m'(x_1, x_2 \otimes x_3) = \\&
 x_1 \otimes x_2 \otimes x_3 +   x_2 \otimes x_1 \otimes x_3+   x_2  \otimes x_3\otimes x_1    - \frac{x_1 x_2}{2\pi i} \otimes  2\pi i \otimes x_3
 -  x_2 \otimes \frac{x_1x_3}{2\pi i} \otimes  2\pi i. \\
&\\
&m'(x_1 \otimes x_2, x_3 \otimes x_4) =\\
&  {x}_1\otimes  {x}_2\otimes {y}_1\otimes {y}_2
+  {x}_1\otimes {y}_1\otimes {x}_2\otimes {y}_2
+  {x}_1\otimes {y}_1\otimes {y}_2\otimes {x}_2\\
&+  {y}_1\otimes {x}_1\otimes {x}_2\otimes {y}_2
+  {y}_1\otimes {x}_1\otimes  {y}_2 \otimes {x}_2
+ {y}_1\otimes {y}_2\otimes {x}_1\otimes {x}_2   \\
&-  {x}_1\otimes  \frac{{x}_2  {y}_1}{2\pi i}\otimes 2\pi i \otimes {y}_2  -
  \frac{{x}_1  {y}_1}{2\pi i}\otimes 2\pi i\otimes  {x}_2\otimes {y}_2 -
   {x}_1\otimes  {y}_1\otimes \frac{{x}_2  {y}_2}{2\pi i} \otimes 2\pi i   \\
&+  \frac{{x}_1 {y}_1}{2\pi i}\otimes 2\pi i \otimes \frac{{x}_2  {y}_2}{2\pi i}\otimes 2\pi i.\\
\end{split}
\ee

We prove that the quasi-shuffle  product $m'$ on ${\rm T}_\bullet(\C)$ is commutative and associative.

Using the projection ${\rm pr}$, see (\ref{PROJJ}), it  induces   quasi-shuffle  product $m$ on ${\cal A}_\bullet$:
\be
\bs
&m(X_1\otimes x_2 \otimes \cdots \otimes x_p, Y_1\otimes y_2 \otimes \cdots \otimes y_q):= \\
&{\rm pr}\circ m'( \log X_1\otimes x_2 \otimes \cdots \otimes x_p,   \log Y_1\otimes y_2 \otimes \cdots \otimes y_q).\\
\end{split}
\ee
We prove that it does not depend on the choice of the logarithms $\log X_1, \log Y_1$.

Another example of the  quasi-shuffle product is discussed in the Appendix B.

\subsection{A map of Hopf algebras $\mathcal{P}: \mathcal{H}_\bullet \rightarrow \mathcal{A}_\bullet$} \la{SECC1.4}

Let us define a map
\be \la{MMAAPPP}
\mathcal{P}: \mathcal{H}_\bullet \rightarrow \mathcal{A}_\bullet.
\ee
We do it in two steps. First,  consider a set  $\widetilde {\cal H}_\bullet$ of {\it split} framed
mixed Hodge-Tate structures. Forgetting the splitting, we get a projection
\be \la{PR}
{\rm p}: \widetilde {\cal H}_\bullet \lra   {\cal H}_\bullet.
\ee
Second, we introduce   our key construction, a map  of sets
$$
\Phi:  \widetilde {\cal H}_\bullet \lra {\rm T}(\C).
$$
We show that the   composition: $\widetilde{\mathcal{P}}:= {\rm pr}\circ \Phi$ descends to a well defined map (\ref{MMAAPPP}). So the map $\mathcal{P}$ is
the unique map making the following diagram commutative:
\be
\begin{gathered}
\xymatrix{
\widetilde {\cal H}_\bullet \ar[r]^-{\Phi} \ar[d]^{\rm p}& {\rm T}(\C) \ar[d]^{{\rm pr}}\\
{\cal H}_\bullet \ar[r]^-{\mathcal{P}}&{\cal A}_\bullet(\C)
}
\end{gathered}
\ee

\paragraph{Construction of the map $\Phi$.}

Let $(f^n; H, s; v_0)$ be a framed mixed Hodge-Tate structure $(f^n, H, v_0)$ with a splitting  $s$.
 Then for any pair
$$
 v_{p} \in \Hom(\mathbb{Q}(p), \gr^W_{-2p}H), ~~~~f^{q} \in \Hom(\gr^W_{-2q}H, \mathbb{Q}(q)), ~~q>p,
 $$ one can define
the period of a framed split Hodge Tate structure $\langle f^q \mid H,s \mid v_p\rangle \in \C$, see    (\ref{period}).

Let us introduce a notation
$$
[ f^q \mid H,s \mid v_p]:=  \frac{\langle f^q \mid H,s \mid v_p\rangle}{(2 \pi i)^{ q-p-1}} \in \C.
 $$

For each  $0<p <  n$, choose a basis $ \{v^\alpha_{p}\}$ in $\Hom(\mathbb{Q}(p), \gr^W_{-2p}H)$. Let $\{f_\alpha^{p}\}$ be the dual basis.

\bd
\la{DECDEF}
The map $\Phi$ is defined  inductively as the unique map satisfying the  identity
\begin{equation}
  \la{REC2a}
  \sum_{0\leq p\leq n}  \sum_{\alpha_p}
[ f^n \mid H,s \mid v^{\alpha_p}_p]\otimes (2 \pi i)^{\otimes(n-p-1)} \otimes \Phi(f_{\alpha_p}^p, H, v_0;s)=0.
\end{equation}
\ed
To make sense out of the formula (\ref{REC2a}), we need to define  its $p=n$ and $p=0$ terms.

1) By definition, the $p=n$ term in the sum (\ref{REC2a})   is just $\Phi(f^n, H, v_0;s)$.

2) By definition, the $p=0$ term in the sum (\ref{REC2a})   is just $ [ f^n \mid H,s \mid v_0] \otimes (2 \pi i)^{\otimes(n-1)}$.

Using 1), we can interpret  the sum (\ref{REC2a}) as an inductive definition of $\Phi( f^n, H,v_0; s)$:
 \be\la{REC}
 \begin{split}
& \Phi(f^n, H,v_0; s):= \\
& -\sum_{0\leq p\leq n-1}  \sum_{ \alpha_p }
  [ f^n \mid H,s \mid v^{\alpha_p}_p]  \otimes (2 \pi i)^{\otimes(n-p-1)} \otimes \Phi(f_{\alpha_p}^p, H, v_0;s).\\
 \end{split}
 \end{equation}

Using the inductive definition, we can write a plain definition of $\Phi( f^n, H,v_0; s)$:
\be \label{MAPPa}
\begin{split}
&\Phi( f^n, H,v_0; s):=\\
&\sum_{\substack{1\leq k\leq n \\ 0=i_0<i_1<\cdots <i_k=n \\  }} \sum_{\alpha_p}
 (-1)^k \bigotimes_{1\leq l \leq k} \left([ f_{\alpha_{i_l}}^{i_l}\mid H,s \mid v^{\alpha_{i_{l-1}}}_{i_{l-1}}] \otimes (2\pi i)^{\otimes (i_l-i_{l-1}-1)}\right)\in {\Bbb C}^{\otimes n}.\\
 \end{split}
 \ee
Here  the sum $\sum_{\alpha_p}$ means the sum  over the basis  vectors
\be
\la{CONV}
 \{f^p_{\alpha_p}, v_p^{\alpha_p}\}, ~~p \in \{0=i_0<i_1<\cdots <i_k=n\}.
\ee

\begin{theorem} \la{Th1.5}  The map $\mathcal{P}:\ \mathcal{H}_\bullet \rightarrow \mathcal{A}_\bullet$ is an isomorphism of graded Hopf algebras.
\end{theorem}

\paragraph{An example: the trilogarithm ${\rm Li}_3(z)$.} Recall the classical polylogarithms:
$$
  {\rm Li}_{n}(x):= \sum_{k=1 }^\infty\frac{x^{k } }{k^{n} }.
$$
The   trilogarithm  ${\rm Li}_3(z)$ defines a rational mixed Hodge-Tate structure that can be represented by the following period matrix:
\[\begin{pmatrix}
  1 \\
  -{\rm Li}_1(z) & 2\pi i\\
  -{\rm Li}_2(z) & 2\pi i \log(z) & (2\pi i)^2\\
  -{\rm Li}_3(z) & (2\pi i) \frac{\log^2(z)}{2} & (2 \pi i)^2 \log (z) & (2\pi i)^3
\end{pmatrix}\]
 Namely, let
 $e_0, e_1, e_2, e_3$ be the natural  basis in the coordinate space $\C^4$. Denote by $C_0, C_1, C_2, C_3$  the vectors in $\C^4$ given by the columns of the matrix,
 counted from the left.
 Then there is the  following  mixed Hodge-Tate structure:
\be
\begin{split}
&H_\mathbb{Q}=\langle C_0, C_1, C_2, C_3\rangle_\mathbb{Q};\\
&W_{-2n} H_{\mathbb{Q}}=\langle C_n, \cdots, C_3\rangle_\mathbb{Q}, ~~W_{-2n +1}= W_{-2n};\\
&H_\mathbb{C}=H_{\mathbb{Q}}\otimes \mathbb{C} = \langle e_0, e_1, e_2, e_3\rangle_\mathbb{C};\\
&F^{-p}=\overline{F}^{-p} =\langle e_0, \cdots, e_{p}\rangle_\mathbb{C}.\\
\end{split}
\ee
Let $v_j:\ \mathbb{Q}(j)\rightarrow H_\mathbb{Q}, \  (2\pi i)^j \mapsto C_j$. Let $f^j$ be the dual of $v_j$.
Then $\langle f^p \mid H \mid v_q\rangle$ is the   entry of the matrix, for example:
$$
\langle f^3 \mid H \mid v_0\rangle=  -{\rm Li}_3(z), ~~ \langle f^2 \mid H \mid v_1\rangle=  \log (z) .
$$

Then
\begin{align*}
  \Phi(f^3, H, v_0)= &- \frac{-{\rm Li}_3(z)}{(2\pi i)^2} \otimes 2\pi i \otimes 2\pi i \\
  &+ \frac{\log^2z}{2\pi i} \otimes 2\pi i \otimes (-{\rm Li}_1(z))\\
  &+ \log z \otimes \frac{-{\rm Li}_2(z)}{2\pi i} \otimes 2\pi i \\
  &- \log z \otimes \log z \otimes (-{\rm Li}_1(z)).
\end{align*}
Therefore we have:
\begin{align*}
  {\cal P}(f^3, H, v_0)= &- \exp\Bigl(\frac{-{\rm Li}_3(z)}{(2\pi i)^2}\Bigr) \otimes 2\pi i \otimes 2\pi i \\
  &+ \exp\Bigl(\frac{\log^2z}{2\pi i}\Bigr) \otimes 2\pi i \otimes (-{\rm Li}_1(z))\\
  &+ z \otimes \frac{-{\rm Li}_2(z)}{2\pi i} \otimes 2\pi i \\
  &- z \otimes \log z \otimes (-{\rm Li}_1(z)).
\end{align*}

  \paragraph{The big period map.}
There is a canonical homomorphism
\be
\begin{split}
&   {\cal A}_n = \C^*\otimes \C^{\otimes n-1} \lra \C^* \otimes \C\\
& A_1 \otimes a_2 \otimes \ldots \otimes a_n \lms (-1) \cdot A_1 \otimes a_2 \ldots a_n.\\
  \end{split}
   \ee
Composing  the   map   ${\cal P}_n: {\cal H}_n\lra {\cal A}_n$ with this map we get a map
$$
P_n: {\cal H}_n \lra \C^* \otimes \C.
$$
By the very definition, this is the big period map defined in \cite{G96}, \cite{G15}. In fact our definition of the map ${\cal P}$ was suggested by the
definition of the   maps $P_n$.

\paragraph{The structure of the paper.} In Section \ref{SEC2.1as} we define and study the Hopf dg-algebra $\mathcal{A}_\bullet(R)$.

  In Section \ref{NS} we study the Hopf dg-algebra $\mathcal{A}_\bullet(R)$ for the Tate dg-algebra $R$ with non-negative grading.
  The crucial example of such a Tate dg-algebra is  the
   de Rham complex $\Omega$ of a manifold $X$ with the Tate line $\Q(1)$.
 We show that the cobar complex of the Hopf dg-algebra $\mathcal{A}_\bullet(\Omega)$
is a dg-model for the rational Deligne cohomology of $X$.

Using this, we suggest in  Section \ref{NS2.3} a dg-model for the category of variations of $\Q$-Hodge-Tate structures on $X$.
There is a  similar  story for a complex regular  algebraic variety $X$.

In Section \ref{SEC2} we study the Hopf dg-algebra $\mathcal{A}_\bullet = \mathcal{A}_\bullet(\C)$ related to the Tate  algebra $\C$ with the Tate line $\Q(1)$.
We define a  map  $\mathcal{P}:\  {\cal H}_\bullet \to {\cal A}_\bullet$.
Then we  prove that it is an morphism of colagebras. This easily implies that the map
$\mathcal{P}$ is an isomorphism.
Then we prove that it commutes with the coproduct, and therefore get    Theorem \ref{Th1.5}.

In Section \ref{SEC4} we prove that $H^0_{\cal D}(\mathcal{A}_\bullet(\Omega))$  the Tannakian   Hopf algebra for the abelian
tensor category of variations of $\Q$-Hodge-Tate structures on $X$.

In Section \ref{Sec5}  we present a variant of the main construction related to the p-adic Hodge theory. 
We start with a version of the notion of the Tate algebra, which we call {\it   Tate $\varphi$-algebra}. 
Given a Tate $\varphi$-algebra $R$ 
we generalize our  construction and define a graded  Hopf algebra ${\cal A}^\varphi_\bullet(R)$. We define a Tate $\varphi$-algebra ${\rm R}_{\rm crys}$ sitting inside of
 Fontaine's period ring $\B_{\rm crys}$. It gives rise to  
a Hopf algebra ${\cal A}^\varphi_\bullet({\rm R}_{\rm crys})$, which  
we relate the    
to   crystalline p-adic mixed Tate Galois representations. 

To treat  the case of  semi-stable mixed Tate Galois representations we outline  in Section \ref{sec5.4} a dg-variant of the story. We introduce the notion of a {\it   Tate $\varphi$-dg algebra} $R^\bullet$, and  
  assign to it a Hopf dg-algebra ${\cal A}^\varphi_\bullet(R^\bullet)$. We interpret Fontaine's ring $\B_{\rm st}$ as a dg-algebra with a differential given by the monodromy $N$, and 
 define a dg-subalgebra ${\rm R}_{\rm st} \subset \B_{\rm st}$. It gives rise   to
 a dg-Hopf algebra ${\cal A}^\varphi_\bullet({\rm R}_{\rm st})$,  which  we 
  relate    
to   semi-stable p-adic mixed Tate Galois representations.  

\paragraph{Acknowledgments.}
This work was
supported by the  NSF grant  DMS-1564385. Its final version was written when the first author enjoyed hospitality of the IHES (Bures sur Yvette) and IAS (Princeton).
 We are grateful to the referee for many useful comments, and to Pierre Colmez for  discussions of the p-adic Hodge theory.

\section{The Hopf dg-algebra $\mathcal{A}_\bullet(R)$ for a   Tate dg-algebra $R$} \label{Hopf_A}

\subsection{The Hopf dg-algebra $\mathcal{A}_\bullet(R)$} \la{SEC2.1as}

\paragraph{Tate dg-algebras.} Let $R$ be a   dg-algebra over a   field $k$ with a differential $d$:\footnote{Later on we assume that $R$ is graded by non-negative integers.}
\be \la{GRAD}
R = \bigoplus_{n=-\infty}^\infty R^n.
\ee
 In particular, $R^0$ is an algebra. We use a shorthand $ |x|$ for the degree ${\rm deg}(x)$ on $R$.

 \bd \la{TATEL} A {\rm Tate line $k(1)$}   in $R$
  is given by the following data:

 \begin{itemize}

 \item

 A 1-dimensional   $k$-vector  subspace $k(1)$ of $d$-constants,  in  the  degree $0$ part of  the center of  $R$:
$$
k(1) \subset {\rm Center}(R)^0, ~~~~dk(1)=0.
$$

\item The subset $k(1)-\{0\}$ belongs to the set $R^*$ of the units of $R$:
\be \la{UNIT}
k(1)-\{0\} \subset R^*.
\ee
\end{itemize}

The data $(R, k(1); d)$ is called a {\rm Tate dg-algebra}.
\ed

Let  $k(0):=k$ and,  using the tensor product over $k$,  set:
$$
k(n):= \otimes^n k(1), ~~~~n> 0.
$$
Thanks to   assumption (\ref{UNIT}),    the set of  inverses  of  non-zero elements    $t \in k(1) -\{0\}$ of the Tate line is
a   $k^*$-torsor. Adding $0$ to it, we get  the {\it inverse Tate line}:
$$
k(-1)\subset R.
$$
Evidently, the product provides an isomorphism $k(-1) \otimes k(1) \lra k$.
Then setting
$$
k(-m):= k(-1)^{\otimes m}, ~~~~m>0,
$$ we get lines $k(n)$ for each integer $n\in \Z$. The   product induces isomorphisms $$
k(n) \otimes k(m) \stackrel{\sim}{\lra} k(n+m).
$$

Given   a {Tate dg-algebra} $(R, k(1))$, we define     a Hopf dg-algebra $\mathcal{A}_\bullet(R)$ with a differential ${\cal D}$.

\paragraph{The   $k$-vector space $\mathcal{A}_\bullet(R)$.}
Denote by      ${\rm T}_\bullet(R)$   the tensor algebra of the $k$-vector space $R$. Set $$
\overline R:=R/k(1).
$$ The  space $\mathcal{A}_\bullet(R)$  is  the quotient of ${\rm T}_\bullet(R)$ by the right ideal generated by the Tate line:
\be
\begin{split}
&{\cal A}_\bullet(R):= {\rm T}_\bullet(R)/ ( k(1)\otimes {\rm T}_\bullet(R)).\\
&{\cal A}_\bullet(R)=  \bigoplus_{n=0}^\infty{\cal A}_n(R), ~~~~
 {\cal A}_0:=k, ~~ {\cal A}_n:=   \overline  R \otimes \underbrace{{R} \otimes
\ldots \otimes {R}}_{\text{ $n-1$ factors }}.\\
\end{split}
\ee

\paragraph{The grading $|\cdot |$.}

Denote by $|\cdot |$ the grading on $R$. We equip $\mathcal{A}_\bullet(R)$   with the {\it   grading}  $|\cdot|$:
\be
\bs
&|x_1 \otimes \ldots \otimes x_n|:=  |x_1|  + \ldots + |x_n|.\\
\end{split}
\ee

\paragraph{The weight grading $w$.} It is induced by the  grading of ${\rm T}_\bullet(R[-1])$:
\be
\bs
&w(x_1 \otimes \ldots \otimes x_n) := n+ |x_1|  + \ldots + |x_n|.
\end{split}
\ee

 \paragraph{The coproduct $\Delta$.} We define a coproduct map on the vector space ${\cal A}_\bullet(R)$ by setting
 $$
\Delta: \overline x_1 \otimes x_2 \otimes \ldots \otimes x_n \lms \sum_{k=0}^n (\overline x_1 \otimes \ldots \otimes x_k) \bigotimes  (\overline x_{k+1} \otimes \ldots \otimes x_n).
$$
Here $\overline x$ is the projection of an element $x\in R$ to $ \overline R $.  The factors with no multiples are equal to $1$.
The   coproduct $\Delta$ is evidently coassociative, and compatible with each of the two gradings.

 \paragraph{The differential ${\cal D}$.}
Using the differential $d$ on $R$,  we define  inductively a map ${\cal D}$ on $\mathcal{A}_\bullet(R)$:
\be \la{DIFFDDa}
{\cal D}(\overline x_1 \otimes \ldots \otimes x_n):= d\overline x_1 \cdot  x_2 \otimes (x_3 \otimes \ldots \otimes x_n) + (-1)^{|x_1|} \overline x_1 \otimes {\cal D}(x_2 \otimes \ldots \otimes x_n).
\ee
It is well defined since $d(k(1))=0$. One checks that  ${\rm deg}({\cal D})=1$, and $w({\cal D})=0$.

 We will  show below that ${\cal D}^2=0$, and that
the coproduct $\Delta$ and  the differential ${\cal D}$ provide   $\mathcal{A}_\bullet(R)$ with a  dg-coalgebra structure for the grading $|\cdot |$.

The definition of the product is more complicated.

 \paragraph{A quasi-shuffle product $m'$ on ${\rm T}_\bullet(R)$.}
Pick a non-zero element
$$
t\in k(1)- \{0\}.
$$

 Using $t$,  we define inductively   a {product} $m'$ on ${\rm T}_\bullet(R)$. It is given by a collection of maps
 $$
 m'_{p,q}:\ {\rm T}_p(R) \otimes {\rm T}_q(R)  \rightarrow {\rm T}_{p+q}(R).
  $$

  First, $1 \in {\rm T}_0(R) = k$ serves as the identity for the product $m'$. This defines
 $m'_{n,0}$  and  $m'_{0,n}$.

Then for $p,q \geq 1$, we define inductively a map   $m_{p,q}$ by setting:
\begin{align*}
  m'(x_1\otimes\cdots \otimes x_p, y_1\otimes\cdots \otimes y_q):=  \quad  &x_1 \otimes m'(x_2\otimes\cdots \otimes x_p, y_1\otimes\cdots \otimes y_q)\\
  &+ (-1)^{\alpha_p\beta_0 }y_1 \otimes m'(x_1\otimes\cdots \otimes x_p, y_2\otimes\cdots \otimes y_q)\\
  &- (-1)^{\alpha_p\beta_0 }x_1 y_1 t^{-1} \otimes t \otimes m'(x_2\otimes\cdots \otimes x_p, y_2\otimes\cdots \otimes y_q).
\end{align*}
It  does not depend on the choice of  $t$: using another $t'=at$, where $a \in k^*$, we get the same $m'$, since the tensor product is over  $k$.
The product is a   degree zero map, thanks to the assumption that the Tate line $k(1)$ is homogeneous. The signs are given by the standard  rules, so:
\be \la{ALB}
\alpha_p := |x_0|+ \ldots + |x_p|, ~~~~\beta_0:= |y_0|.
\ee
 We will show that  $m'$
 defines an associative product on ${\rm T}_\bullet(R)$.

\paragraph{The   product $m$.} Although the product $m'$ and coproduct $\Delta'$ in ${\rm T}_\bullet(R)$ are \underline{non-compatible}
since  they do not satisfy the Hopf axiom,
they induce a  Hopf algebra structure on $\mathcal{A}_\bullet(R)$.

Precisely,  we define a product $m$ on $\mathcal{A}_\bullet(R)$ by the following inductive formula:

\be \label{EQPin}
  \begin{split}
  m(\overline x_1\otimes x_2 \otimes \cdots \otimes x_p, \overline y_1\otimes y_2 \otimes \cdots \otimes y_q)&:= \overline x_1 \otimes m'(  x_2\otimes\cdots \otimes x_p, y_1\otimes\cdots \otimes y_q)\\
  &+ (-1)^{\alpha_p\beta_0}\overline y_1 \otimes m'(  x_1\otimes\cdots \otimes x_p,  y_2\otimes\cdots \otimes y_q)\\
  &- (-1)^{\alpha_p\beta_0}\overline {x_1 y_1t^{-1}} \otimes t\otimes m'(  x_2\otimes\cdots \otimes x_p,  y_2\otimes\cdots \otimes y_q).
\end{split}
\ee
 Notice that we use the product $m'$ on the right.  Here
  $x_1, y_1$ are arbitrary    lifts of the elements $\overline x_1, \overline y_1 \in \overline {R^0}$ to $R^0$.
A key point is that, as we  will prove   below, the product  does not depend on the choices of the lifts.
For example, the simplest non-trivial component of $m$ is:
$$
m(\overline x_1 , \overline y_1):= \overline x_1 \otimes  y_1+ \overline y_1 \otimes  x_1 - \overline{x_1y_1t^{-1}}\otimes t.
$$

\paragraph{The quasi-shuffle formula for the product $m'$.}
Let us give a formula for the product $m'$ of $n$ elements:
\be \la{QSS1}
m'(x^{(1)}_1\otimes\cdots \otimes x^{(1)}_{p_1}, \ldots , x^{(n)}_1\otimes\cdots \otimes x^{(n)}_{p_n}).
\ee

Let us define a quasi-shuffle   of the ordered sets
\be \la{QSS}
\{s^{(1)}_1,\cdots ,s^{(1)}_{p_1}\}, \ldots , \{s^{(n)}_1,\cdots , s^{(n)}_{p_n}\}.
\ee
Consider   a configuration of  distinct points $z_1, ..., z_N$ on the real line, such that each point of the configuration is labeled by a non-empty subset   of the  union of the sets (\ref{QSS}),
so that

1. The labels of each point $z$ contain no more then one element of each of the sets.

2. For each $1 \leq k \leq n$, the line order of the points whose labels contain the elements $s^{(k)}_1,\cdots , s^{(k)}_{p_k}$   coincides with the  order of the  set
$s^{(k)}_1,\cdots , s^{(k)}_{p_k}$.

Such a labeled configuration of points, considered modulo isotopy,  encodes a quasi-shuffle of the $n$ ordered sets  (\ref{QSS}).

 To define the product (\ref{QSS1}), we assign to  each quasi-shuffle  an element of $R^{\otimes (p_1+\ldots + p_n)}$.
 We assign to each point $z$ labeled by $a$ elements $s_1, ..., s_a$, assigned to the elements $x_1, ..., x_a$ of $R$ by matching the data in (\ref{QSS1}) and (\ref{QSS}), the element
 $$
 (-1)^{a-1}\cdot x_1 \ldots x_a t^{-(a-1)} \otimes t^{\otimes (a-1)} \in R^{\otimes a}.
 $$
Then we multiply these elements  following the line order of the points $z_j$.  The quasi-shuffle product (\ref{QSS1}) is defined as the sum of the obtained elements over all quasi-shuffles of (\ref{QSS}).

To show that the quasi-shuffle product $m'$ is well defined one needs to show that the product of two such products is again given by the quasi-shuffle formula.
This boils down to checking
\begin{align*}
 (-1)^{a-1}m'\Bigl( x_1 \ldots x_a t^{-(a-1)} \otimes t^{\otimes (a-1)}, ~   y\Bigr) &=  (-1)^{a}x_1 \ldots x_a yt^{-a} \otimes t^{\otimes a}\\
 &+ (-1)^{a-1} (-1)^{\alpha_a\beta}y \otimes x_1 \ldots x_a t^{-(a-1)} \otimes t^{\otimes (a-1)}\\
 &+ (-1)^{a-1}  x_1 \ldots x_a t^{-(a-1)} \otimes t^{\otimes (a-1)} \otimes y\\
 &\in R^{\otimes (a+1)}.
\end{align*}
Here $(-1)^{\alpha_a\beta}$ is the standard sign obtained by moving $y$ through $x_1...x_p$. Let us prove this.  Colliding the point $y$ with the point labeled by $x_1 \ldots x_a t^{-(a-1)}$, and using the fact that $t$ is in the center, we get the first term on the right hand side. Putting $y$ on the very left or on the very right, we get the second and third term.
 The other terms of the quasi-shuffle product vanish. To see this, notice  that the   collision of $y$ with $(j+1)$-st factor $t$ delivers a factor
  $$
  (-1)^{a}x_1 \ldots x_a t^{-a} \otimes \underbrace{t\otimes \ldots \otimes t}_{\mbox{$j$ times}} \otimes \underbrace{y \otimes t} \otimes \ldots
  $$
  Here the second underbraced factor $y \otimes t$ appears as $-tyt^{-1} \otimes t$.

  On the other hand moving $y$ a bit to the left of the $(j+1)$-st factor $t$ gives
 $$
  (-1)^{a-1}x_1 \ldots x_a t^{-a} \otimes \underbrace{t\otimes \ldots \otimes t}_{\mbox{$j$ times}} \otimes y \otimes t\ldots
  $$
  So the two terms cancel. This proves, by the induction, that the product $m'$ is well defined. Then it is obvious that it is associative, and commutative if $R$ is commutative.

  A formal arguments using the inductive definition will be also presented below.
\bt \la{TH2.2}
The vector space  ${\cal A}_\bullet(R)$ has a  Hopf dg-algebra structure with the deconcatenation coproduct  $\Delta$, the quasi-shuffle product $m$,    
 the  degree ${\rm deg}$, and the differential ${\cal D}$, with  ${\rm deg}({\cal D})=1.$ 
 
 The Hopf dg-algebra  ${\cal A}_\bullet(R)$ has  an extra weight grading $w$. 
 One has $w({\cal D})=0$.
  \et

\begin{proof} Let us start the proof of Theorem \ref{TH2.2}.

\paragraph{Properties of the quasi-shuffle products on ${\rm T}_\bullet(R)$ and $\mathcal{A}_\bullet(R)$.}

The following Lemma describes the role of the Tate line $k(1)\subset R$  in the multiplication.

\bl \la{L4.2}
$m'_{1+p, q}(t\otimes x_1 \otimes \cdots \otimes x_p, y_1\otimes\cdots y_q)=t\otimes m'_{p,q}(x_1\otimes \cdots \otimes x_p, y_1\otimes\cdots \otimes y_q).$
\el

\begin{proof}
  For $q=0$, it is trivial.
  For $q>0$,
  \begin{align*}
    m'_{1+p, q}(t\otimes x_1 \otimes \cdots \otimes x_p, y_1\otimes\cdots y_q)
    =&\quad t\otimes m'_{p,q}(x_1\otimes \cdots \otimes x_p,y_1\otimes\cdots \otimes y_q)\\
    &+ (-1)^{\alpha_p\beta_0 }y_1\otimes m'_{1+p, q-1}(t\otimes x_1 \otimes \cdots \otimes x_p, y_2\otimes\cdots y_q)\\
    &-(-1)^{\alpha_p\beta_0 }y_1\otimes t \otimes m'_{p, q-1}(x_1 \otimes \cdots \otimes x_p, y_2\otimes\cdots y_q).
  \end{align*}

  By the inductive assumption, the last two terms cancel out.
\end{proof}

\bc \la{C4.3}
  The following equations hold:

  (1) $m'_{p,q}(t^{\otimes p}, t^{\otimes q})=t^{\otimes (p+q)}$;

  (2) $m'_{p_1+p_2, q}(t^{\otimes p_1}\otimes x_1 \otimes\cdots \otimes x_{p_2}, y_1\otimes \cdots \otimes y_q)=t^{\otimes p_1} \otimes m'_{p_2, q_1}(x_1\otimes\cdots \otimes x_{p_2}, y_1\otimes\cdots \otimes y_q)$.
\ec

Using the part (2) in Corollary \ref{C4.3} we have the following lemma, saying that the product can be calculated in ``blocks'' of the form $x\otimes t^{\otimes p}$.

\bl \la{L4.4}  One has, setting  $
\alpha_p :=  |x_0|+ \ldots + |x_p|, ~~~~\beta_0:= |y_0|$:
  \begin{align*}
    &m'_{p_1+p_2, q_1+q_2}(x_0\otimes t^{\otimes p_1-1}\otimes x_1\otimes\cdots \otimes x_{p_2}, y_0\otimes t^{\otimes q_1-1}\otimes y_1\otimes\cdots y_{q_2})\\
    =&\quad x_0\otimes t^{\otimes p_1-1} \otimes m'_{p_2, q_1+q_2}(x_1\otimes\cdots \otimes x_{p_2}, y_0\otimes t^{\otimes p_1-1}\otimes y_1\otimes\cdots y_{q_2})\\
    &+ (-1)^{\alpha_{p_2}\beta_0 }y_0\otimes t^{\otimes q_1-1} \otimes m'_{p_1+p_2, q_2}(x_0\otimes t^{\otimes p_1-1}\otimes x_1\otimes\cdots \otimes x_{p_2}, y_1\otimes\cdots y_{q_2})\\
    &- (-1)^{\alpha_{p_2}\beta_0 }x_0 y_0 t^{ -(p_1+q_1-1)}\otimes t^{\otimes p_1+q_1-1}\otimes m'_{p_2, q_2}(x_1\otimes\cdots \otimes x_{p_2},y_1\otimes\cdots y_{q_2}).
  \end{align*}
\el

Lemma \ref{L4.4} easily implies that the product $m'$ is associative.

 There are natural surjective projections   $\ R^{\otimes n}\rightarrow \mathcal{A}_n(R)$. Consider their sum:
$$
{\rm pr}:\ {\rm T}_\bullet(R) \longrightarrow \mathcal{A}_\bullet(R).
$$
The  product $m$ on $\mathcal{A}_\bullet(R)$ is induced by the one $m'$ on ${\rm T}_\bullet(R)$, so that we have
a commutative diagram:

\centerline{
\xymatrix{
{\rm T}_\bullet(R)\otimes {\rm T}_\bullet(R) \ar[r]^{~~~~~~~m'} \ar[d]_{{\rm pr} \otimes {\rm pr}}& {\rm T}_\bullet(R) \ar[d]_{  {\rm pr}}  \\
\mathcal{A}_\bullet(R)\otimes \mathcal{A}_\bullet(R)  \ar[r]^{~~~~~~~m}  & \mathcal{A}_\bullet(R)
}
}

The product $m$ is well-defined by Lemma \ref{L4.2}.
\bp
  $(\mathcal{A}_\bullet(R),m)$ is an associative graded algebra. It is commutative if and only if the algebra $R$ is commutative.
\ep

\begin{proof} Since the projection ${\rm pr}: {\rm T}_\bullet(R) \lra {\cal A}_\bullet(R)$ is surjective, this follows from the   properties of the product $m'$, which has been already established. \end{proof}

\paragraph{The Hopf dg-algebra structure on ${\cal A}_\bullet(R)$.}

Let $B$ be a  $k$-vector space graded by non-negative integers, with $B_0=k$.   The product $m$ and  coproduct $\Delta$  form a  dg-bialgebra if and only if
  the Hopf axiom holds:
\be \la{HA}
\Delta(X\ast Y) - \Delta(X) \ast_2   \Delta(Y) =0.
\ee
Here we use $\ast$ for the product $m$, and $\ast_2$ for the induced product on $B \otimes B$.
Precisely, the Hopf axiom tells that the following diagram must be commutative, where    $\tau(x\otimes y)= (-1)^{|x||y|}y\otimes x$:

\vskip 2mm
~~~~~~~~~~~~~~~~~~~~~~~~\xymatrix{
B \otimes B \ar[r]^{m} \ar[d]_{\Delta \otimes \Delta}& B \ar[r]^{\Delta} & B \otimes B\\
B\otimes B \otimes B \otimes B \ar[rr]^{\Id \otimes \tau \otimes \Id} & & B \otimes B \otimes B \otimes B \ar[u]_{m\otimes m}
}

\paragraph{Warning.}  $({\rm T}_\bullet(R), m^\prime, \Delta)$ is not a bialgebra since already for $x,y \in R_0 \subset {\rm T}_1(R)$ we have:
\be
\begin{split}
&x\ast y = x\otimes y + y \otimes x - xyt^{-1} \otimes t;\\
&\Delta(x\ast y) - \Delta(x) \ast_2   \Delta(y) = - xyt^{-1} \otimes t. \\
\end{split}
\ee

 Lemma \ref{L4.8}  claims that the difference (\ref{HA}) lies in a subspace ${\rm T}_\bullet(R)\otimes {\rm Ker} ({\rm pr})$ of ${\rm T}_\bullet(R)\otimes {\rm T}_\bullet(R)$.

\bl \la{L4.8} We have:
 \be
 \Big((\Delta\circ m')-(m'\otimes m')\circ \big((\Id \otimes \tau \otimes \Id)\circ (\Delta \otimes \Delta)\big)\Big)\big({\rm T}_\bullet(R)\otimes {\rm T}_\bullet(R)\big)\subset {\rm T}_\bullet(R)\otimes {\rm Ker} ({\rm pr}).
 \ee
\el

\begin{proof} The shuffle product and the deconcatenation coproduct on the tensor algebra   of a vector space  satisfy the Hopf axiom.
In our case   the quasi-shuffle product produces a factor $x_iy_j t^{-1}\otimes t$ for each "colliding pair" $(x_i, y_j)$.
Applying the deconcatenation  between
$x_iy_jt^{-1}$ and $t$ we get an extra term $( ...\otimes  x_iy_j t^{-1}) \bigotimes (t \otimes ...)$.  Since the second factor starts with $t$, we get the claim.
\end{proof}

Since ${\rm T}_\bullet(R)\otimes {\rm Ker} ({\rm pr})\subset {\rm Ker} ({\rm pr}\otimes  {\rm pr})$,   Lemma \ref{L4.8}  implies that $\mathcal{A}_\bullet(R)$ satisfies the Hopf axiom. 
So $(\mathcal{A}_\bullet(R),m,\Delta)$ is a dg-bialgebra.
Being $\Z_{\geq 0}$-graded and connected,  it has a unique Hopf algebra structure.

\paragraph{Properties of the differential ${\cal D}$.}

  Let us show that ${\cal D}^2=0$.
    Observe that
\be
\bs
 &{\cal D} \circ {\cal D}(f \otimes g\otimes h)  =   {\cal D} (df \cdot g\otimes h + (-1)^{|f|}f\otimes dg\cdot h)  =\\
 & (-1)^{|f|+1}df \cdot dg \cdot h + (-1)^{|f|}df \cdot dg \cdot h=0.
   \end{split}
   \ee
  There is a similar argument when ${\cal D} \circ {\cal D}$ affects disjoint pairs $f_1 \otimes f_2 \otimes \ldots \otimes g_1 \otimes g_2$ in the tensor product.

  Let us show that  ${\cal D}$   satisfies the Leinbiz rule for the product
  on ${\rm T}_\bullet(R)$. Using   $dt=0$,
  $$
  {\cal D}(m'(f, g)) =     {\cal D}(f \otimes g + (-1)^{|f|\cdot |g|}g \otimes f - f gt^{-1} \otimes t) = (df) g + (-1)^{|f|\cdot |g|}(dg)f - d(fg) =0.
  $$
  The argument in general is reduced to this one. The claim for ${\cal A}_\bullet(R)$ follows from this.

Let us show that   ${\cal D}$ is compatible
   with the coproduct. Take an element $f_1 \otimes \ldots \otimes f_n$. Applying  ${\cal D}$, we pick a pair of consecutive factors $f_i \otimes f_{i+1}$ and map  it to
   $df_i \cdot f_{i+1}$.   Applying then the coproduct we cut between $f_k$ and $f_{k+1}$ where either $k<i$ or $k>i$. On the other hand, applying the coproduct first we cut
   between $f_k$ and $f_{k+1}$. But if $k=i$, we do not get a term with  $df_i \cdot f_{i+1}$.

  .
  Theorem \ref{TH2.2} is proved.\end{proof}

\bl
The graded Hopf algebra  $H^0_{\cal D}({\cal A}_\bullet(R))$
 is  realized as the subspace of ${\cal D}$-cycles in  ${\cal A}_\bullet({R^0})$:
$$
H^0_{\cal D}({\cal A}_\bullet(R)) = \{x \in {\cal A}_\bullet(R^0)~|~ {\cal D}(x)=0 \} \subset  {\cal A}_\bullet(R^0).
$$
 \el

 \begin{proof}  Since ${\cal A}_\bullet(R)$ is a Hopf dg-algebra,  $H^0_{\cal D}({\cal A}_\bullet(R))$ is a graded Hopf algebra.
   The degree $0$ part of ${\cal A}_\bullet(R)$ is   ${\cal A}_\bullet(R^0)$, and the degree $-1$ part is zero.
   \end{proof}

\paragraph{Remark.} The tensor algebra ${\rm T}_\bullet(R)$ has a standard commutative Hopf dg-algebra structure: its differential ${\rm D}$ is
induced on the tensor product by the differential $d$ on $R$.
The product is the shuffle product. The coproduct is given by the deconcatenation. However it does not induce a Hopf dg-algebra structure
on the quotient  ${\cal A}_\bullet(R)$.
Furthermore, the differential ${\cal D}$ on ${\cal A}_\bullet(R)$ requires the algebra structure on $R$, while ${\rm D}$ does not.

\paragraph{The cobar complex.}   

Given a dg-coalgebra $({\rm A}, \Delta, d)$, there is the   cobar complex:
$$
{\rm Cobar}_\bullet({\rm A}):=~~ {\rm T}_\bullet({\rm A}[-1]).
$$
Its differential is a sum of two commuting differentials: the one is the restricted coproduct $\Delta'$; the other is induced by  $d$, via  the Lebniz rule on the   algebra ${\rm T}_\bullet({\rm A}[-1])$.
If the algebra ${\rm A}$ is commutative, the shuffle product  makes the cobar complex into   a commutative dg-algebra.

 The   cobar complex   ${\rm Cobar}_\bullet({\cal A}_\bullet(R))$ of the   Hopf dg-algebra ${\cal A}_\bullet(R)$
 is graded by the weight.
   Denote by ${\rm Cobar}_n({\cal A}_\bullet(R))$ its weight $n$ part.

   So   a commutative Tate dg-algebra $R$ provides a   commutative dg-algebra,   with an additional grading by the weight preserved by the differential:
\be \la{MCOz}
{\rm Cobar}_\bullet({\cal A}_\bullet(R)) = \bigoplus_{n=0}^\infty{\rm Cobar}_n({\cal A}_\bullet(R)).
\ee

 \bl \la{EXL1a} Assume that $R^0=k$. Then
the complex ${\rm Cobar}_n({\cal A}_\bullet(R))$ is isomorphic to
the complex
$$
\overline R\bigotimes \Bigl(R \stackrel{}{\lra} \overline R \Bigr)^{\otimes n-1}.
$$
  The complex ${\rm Cobar}_n({\cal A}_\bullet(R))$ is a resolution of $R/k(n)$.\el

\begin{proof}
The first claim  follows from the very definitions. The second  follows from the first. \end{proof}

 In Section \ref{NS} we explore the analog of Lemma \ref{EXL1a} for Tate dg-algebras   graded  by non-negative integers.

 \subsection{Commutative Tate  dg-algebras and Deligne complexes} \la{NS}
 Assume now that $(R, d)$ is a  commutative dg-algebra, graded by non-negative integers:
 $$
 R = R^0 \oplus R^1 \oplus R^2 \oplus \ldots .
 $$
  Then a Tate algebra structure on $R$ gives rise to
  an analog of the  weight $n$ Deligne  complex:
\be \la{WTNDC}
 \xymatrix{\underline k^\bullet_{R; {\cal D}}(n)  :=&  k(n)  \ar[r]^{} &R^0\ar[r]^{d}& R^1 \ar[r]^{d}& \ldots \ar[r]^{d}& R^{n-1}. & }
\ee
Here $k(n)$ is in the degree zero. One can replace the first two terms by $\overline {R^0}(n-1)$.

\paragraph{Example.}
The   weight $n$ Deligne complex
  on a complex manifold $X$ is a complex of sheaves
$$
\begin{array}{ccccccccccccccc}
{ \underline \Q}^\bullet_{{\cal D}}(n):=~~~~&{  \Q}(n)&\lra &\Omega^0& \stackrel{d}{\lra}&\Omega^1& \stackrel{d}{\lra}&\ldots&\stackrel{d}{\lra} &\Omega^{n-1}.&
\end{array}
$$
The rational weight $n$ Deligne cohomology of $X$ are the cohomology of $X$ with coefficients in this complex of sheaves.
Therefore when $R =(\Omega, d)$ is the de Rham complex of sheaves  on a complex manifold $X$,  and $k(1) = \Q(1)$,  the complex (\ref{WTNDC})
coincides with the  rational weight $n$ Deligne complex of sheaves on   $X$.

  \bt \la{TH2.11} The complex   ${\rm Cobar}_n(\mathcal{A}_\bullet(R))$
     is quasi-isomorphic to the weight $n$ Deligne complex (\ref{WTNDC}).
  \et

 \begin{proof}  We will use the traditional notation $~|~$ for the tensor product in the cobar complex,   keeping the sign $\otimes$ for the tensor product everywhere else,
  including ${\cal A}_\bullet(R)$.

A   term in  cobar complex (\ref{MCOz}) is given by a tensor product  of $\overline {R^0}$, $R^0$, and
 $R^p$ where $p>0$.

  Let us consider a decreasing filtration ${\cal F}^\bullet$ on the cobar complex such that   
 $ {\cal F}^p$ is a  sum of the terms where the total degree of the factors  $R^s$,  $s>0$, entering  the term is   $\geq p$. We calculate the cohomology of the cobar complex using the spectral sequence
 for the filtration ${\cal F}^\bullet$.

 \bp \la{P5.2}
 There is a quasi-isomorphism
 $$
R^p(n-p-1)  \stackrel{\rm quis}{\lra} {\rm gr}_{\cal F}^p \Bigl({\rm Cobar}_\bullet({\cal A}_\bullet(R))\Bigr).
 $$
 The induced differential in the spectral sequence coincides with the de Rham differential
 $$
d:  R^p(n-p-1) \lra  R^{p+1}(n-p-2).  $$
 \ep

 \begin{proof}

Recall   the reduced  graded   algebra $\overline R$:
    $$
  \overline{R} = \overline {R^0} \oplus R^1 \oplus \cdots.
  $$

  The differential in the cobar complex ${\rm Cobar}_\bullet({\cal A}_\bullet(R))$ has two components: the one induced by the differential ${\cal D}$, and the one
 provided by the restricted coproduct in ${\cal A}_\bullet(R)$.
 Passing to  ${\rm gr}_{\cal F}^p$, the  differential ${\cal D}$ vanishes.  We denote the weight $n$ part of a weight graded space $A$ by $A_{w=n}$. Then there is a canonical isomorphism of complexes, generalizing Lemma \ref{EXL}:

  $$
  \gr_\mathcal{F}^\bullet ({\rm Cobar}_n(\mathcal{A}_\bullet(R)))= \gr_\mathcal{F}^\bullet
  \Bigl(\bigoplus_{i=0}^\infty \overline{R} \otimes (R\longrightarrow \overline{R})^{\otimes i}\Bigr)_{w=n}.
  $$
  Indeed, each factor $\overline R$ coming from $(R\longrightarrow \overline{R})^{\otimes i}$ matches one of the  bars $|$  in the cobar complex, in agreement
  with the fact that both  this factor and $|$  have the degree $1$.

  Since the complex $R \longrightarrow \overline{R}$ is a quasi-isomorphic to $k(1)$, we have:
    $$
  \overline{R} \otimes (R  \longrightarrow \overline{R})^{\otimes i} ~~\overset{\rm quis}\sim ~~\overline R(i).
  $$
  So
  $$
  \gr_\mathcal{F}^p  \Bigl(\bigoplus_{i=0}^\infty \overline{R} \otimes (R \longrightarrow \overline{R})^{\otimes i}\Bigr)_{w=n}~~ \overset{\rm quis}\sim ~~
  \overline R^p(n-p-1).
  $$
  The spectral sequence differential   $ \overline R^p(n-p-1) \lra \overline R^{p+1}(n-p-2)$ is given, for $p>0$, by:
  \be
  \bs
  & {\cal D}: R^p \otimes  k(1)^{\otimes n-p-1}  \longrightarrow R^{p+1} \otimes k(1)^{\otimes n-p-2}, \\
  &\omega_p \otimes t^{\otimes n-p-1} \lms t  d \omega_p  \otimes t^{\otimes n-p-2}. \\
  \end{split}
  \ee
So we recover the Deligne complex  $\underline k^\bullet_{R; {\cal D}}(n)$ for the Tate dg-algebra $R $, with the differential  $td$:
 $$
 \Bigl(\bigoplus_{i=0}^\infty \overline{R } \otimes (R  \longrightarrow \overline{R })^{\otimes i}\Bigr)_{w=n}~~ \overset{\rm quis}\sim ~~ \overline {R^0}(n-1) \lra R^1(n-2) \lra \ldots \lra
 R^{n- 1}.  $$
 Proposition \ref{P5.2}, and therefore  Theorem \ref{TH5.1} are proved. \end{proof}  \end{proof}

\subsection{A commutative dg-model for the rational Deligne cohomology}  \la{NS2.3}

Let $(\Omega, d)$ be the sheaf of commutative dg-algebras of holomorphic forms on a complex manifold $X$. Then
${\cal A}_\bullet(\Omega )$ is a sheaf of commutative Hopf dg-algebras on $X$, with an extra   weight grading.
Its cobar complex   is a sheaf of commutative dg-algebras on $X$,  with an extra   weight grading:
\be \la{MCOzaaa}
{\rm Cobar}_\bullet({\cal A}_\bullet(\Omega)) = \bigoplus_{n=0}^\infty{\rm Cobar}_n({\cal A}_\bullet(\Omega)).
\ee

When $X$ is a regular complex algebraic variety, we take its compactification with normal crossing divisor at infinity, and consider a similar complex where the de Rham complex $\Omega$ is replaced by the  de Rham complex $\Omega_{\log}$ of forms
 with logarithmic singularities at infinity. In this case we consider the Hopf dg-algebra ${\cal A}_\bullet(\Omega_{\log}) $.

\bt \la{TH5.1}
 Let $X$ be  a complex manifold. Then
 the weight $n$ cobar complex  ${\rm Cobar}_n({\cal A}_\bullet(\Omega)) $  is quasi-isomorphic to  the   Deligne   complex of sheaves $\underline \Q^\bullet_{\cal D}(n)$ of $X$.

   So the  commutative dg-algebra (\ref{MCOzaaa}) is  a dg-model for the rational Deligne cohomology of $X$.

   Similar results hold for a regular complex variety  and its  Hopf dg-algebra ${\rm Cobar}_\bullet({\cal A}_\bullet(\Omega_{\log}))$.
 \et

 \begin{proof}  This is a special case of Theorem \ref{TH2.11}. \end{proof}

 \paragraph{Examples.}

1.  The weight 2 part of the cobar complex ${\rm Cobar}_2({\cal A}_\bullet(\Omega ))$ looks as follows:
 \be \la{MCO}
\begin{gathered}
\xymatrix{
    {\cal O}^* \otimes {\cal O} \ar[d]^{\cal D}\ar[r]&  {\cal O}^* ~|~  {\cal O}^*   \\
   {\Omega}^1 &      }
\end{gathered}
\ee
 Replacing   the top row by the quasi-isomorphic ${\cal O}^*(1)$, we get the complex $\underline \R_{\cal D}(2)= {\cal O}^*(1) \lra \Omega^1$.

2. The weight 3 part of the cobar complex ${\rm Cobar}_3({\cal A}_\bullet(\Omega ))$ looks as follows:
\be \la{MCO3}
\begin{gathered}
\xymatrix{
   {\cal O}^* \otimes {\cal O} \otimes {\cal O}\ar[d]^{\cal D}\ar[r]&  {\cal O}^* ~|~
{\cal O}^*\otimes {\cal O}  \bigoplus {\cal O}^*\otimes  {\cal O}~|~   {\cal O}^*\ar[d]^{\cal D} \ar[r] &{\cal O}^*~|~   {\cal O}^*~|~  {\cal O}^* \\
    {\Omega}^1 \otimes {\cal O}  \bigoplus {\cal O}^*   \otimes {\Omega}^1\ar[r]^{} \ar[d]^{\cal D}&    {\cal O}^* ~|~ {\Omega}^1 \bigoplus   {\Omega}^1  ~|~  {\cal O}^*   \\
  \Omega^2}
\end{gathered}
\ee

Replacing each row by  quasi-isomorphic complexes, we get the weight $3$ Deligne  complex:

\be \la{MCO1}
\begin{gathered}
\xymatrix{
{\cal O}^*(2)  \ar[r]^{}  &
  {\Omega}^1(1)  \ar[r]^{}   &
   \Omega^2}.
\end{gathered}
\ee

\bd The dg-category of complexes of variations of Hodge-Tate structures on a complex manifold $X$
is   the dg-category of dg-comodules over the Hopf dg-algebra ${\cal A}_\bullet(\Omega )$.
\ed

This definition is discussed in the next Section.

 \subsection{Tate dg-algebras and mixed Tate categories } \la{SECC1.1}

The category of comodules over a  Hopf algebra $A$ is an abelian tensor category.
 The same is true if $A$ is a Hopf algebra in a tensor category ${\cal C}$.
 In the examples below the category ${\cal C}$ is either the category of graded vector spaces, or the category of complexes.

 Therefore Theorem \ref{TH2.2} tells us that a Tate dg-algebra $R$ gives rise to a tensor dg-category
 ${\cal M}^\bullet(R)$  of dg-comodules over the Hopf dg-algebra $\mathcal{A}_\bullet(R)$.

 Denote by  ${\cal M}(R)$  the abelian tensor category    of graded $H^0_{\cal D}(\mathcal{A}_\bullet(R))$-comodules.

  The cohomology of  dg-comodules over  $\mathcal{A}_\bullet(R)$ are $H^0_{\cal D}(\mathcal{A}_\bullet(R))$-comodules,
 providing   functors
 $$
 H^i: {\cal M}^\bullet(R)\lra  {\cal M}(R).
 $$
   If $R$ is just an algebra, that is the grading on $R$ is trivial, $R=R^0$, then $H^0_{\cal D}(\mathcal{A}_\bullet(R)) = \mathcal{A}_\bullet(R)$.

 \paragraph{Mixed Tate categories arising from  commutative Tate  algebras.}
  Assume that $R$ is a commutative algebra. Then the tensor   category   ${\cal M}(R)$ is a {\it mixed Tate category}, see Appendix \ref{SSEc5}.
 Recall that a mixed Tate category ${\cal T}$ is an abelian Tannakian $k$-category with invertible object $k(1)_{\cal T}$ such that the tensor powers of $k(1)_{\cal T}$ and its dual $k(-1)_{\cal T}$
provide   non-isomorphic simple objects $k(n)_{\cal T}$, $n \in \Z$. Any simple object of the category is isomorphic to one of them.

  The category ${\cal M}(R)$ assigned to a mixed Tate algebra $(R, k(1))$ has in addition to this the following two  properties:

\begin{enumerate}

\item  It has the homological dimension $1$: 
\be \la{CONDI1}
\bs
&{\rm Ext}^k_{{\cal M}(R)}(k(0)_{\cal T}, k(n)_{\cal T}) =0, ~~~~\forall k>1. \\
\end{split}
\ee

\item Its ${\rm Ext}^1$'s are determined by
\be \la{CONDI2}
\bs
&{\rm Ext}^1_{{\cal M}(R)}(k(0)_{\cal T}, k(n)_{\cal T}) =  \left\{ \begin{array}{lll}
R/k(n) & \mbox{ if } &  n > 0,\\
0
& \mbox{ if } &  n \leq 0\\\end{array}\right. \\
\end{split}
\ee
In particular,  the spaces (\ref{CONDI2}) are canonically isomorphic as $k$-vector spaces.
\end{enumerate}
This follows   from Theorem \ref{TH2.11}, since  if $R=R^0$, the complex (\ref{WTNDC}) is quasi-isomorphic to $R/k(n)$.

\bl A mixed Tate category satisfying the conditions (\ref{CONDI1})-(\ref{CONDI2}) where $(R, k(1))$ is a commutative Tate algebra,
 is  equivalent to the  category ${{\cal M}(R)}$.

\el

\begin{proof} Both categories are determined by the structure of the Ext groups.
\end{proof}

A much more interesting result,  proved in Section \ref{SEC2}, is that the category of mixed $\Q$-Hodge-Tate structures, which is a
 mixed Tate category satisfying  conditions (\ref{CONDI1})-(\ref{CONDI2}) for the Tate algebra  $(\C, \Q(1))$, is {\bf canonically} equivalent
 to the category  ${{\cal M}(\C)}$.
 This  recovers Example 1) in Section \ref{SECC1.1a}.

\section{The period morphism  ${\cal P}: {\cal H}_\bullet \stackrel{}{\longrightarrow}  {\cal A}_\bullet$} \la{SEC2}

\subsection{The period map $\cal P$}

 \paragraph{The period operator.}

Let $H$ be a mixed Hodge-Tate structure over $\mathbb{Q}$. There is an isomorphism:
\[H_\mathbb{C}=\bigoplus_p F^p H_\mathbb{C}\cap W_{2p}H_\mathbb{C}.\]
Furthermore, the following canonical map is an isomorphism:
\[F^p H_\mathbb{C}\cap W_{2p}H_\mathbb{C} \overset{\sim}{\rightarrow} \gr_{2p}^W H_\mathbb{Q} \otimes_\mathbb{Q} \mathbb{C}.\]
Using the isomorphisms above we get a canonical isomorphism:
\[S_{HT}:\ H_\mathbb{C} \overset{\sim}\rightarrow \bigoplus_p \gr_{2p}^W H_\mathbb{C}.\]
On the other hand a splitting of the weight filtration on $H_\mathbb{Q}$ provides an isomorphism
\[S_W:\ \bigoplus_p \gr_{2p}^W H_\mathbb{C} \overset{\sim}\rightarrow H_\mathbb{C}.\]
Then the composition of these isomorphisms provides a map
\[S_{HT} \circ S_W \in \End(\bigoplus_p \gr_{2p}^W H_\mathbb{C}),\]
called the period operator. It is determined by the splitting.

For any split mixed Hodge-Tate structure $(H,s)$ and an $(i,j)$-framing $$
v \in \Hom(\mathbb{Q}(i), \gr^W_{-2i}H), ~f\in \Hom(\gr^W_{-2j}H, \mathbb{Q}(j))
$$ we define the period
\be
\label{period}
 \langle f\mid H, s \mid v\rangle := f\circ S_{HT}\circ S_W \circ v \in \Hom_\mathbb{C}(\mathbb{Q}(i)_\mathbb{C}, \mathbb{Q}(j)_{\mathbb{C}})\cong \mathbb{C}.
\ee
The last isomorphism is normalized so that it sends the map $  (2\pi i)^i \mapsto  (2\pi i)^n$ to $(2\pi i)^{n-i}\in \mathbb{C}$.
 For a given $(H,s)$, this map is linear in $f\otimes v$.
We also use the notation
$$
[f\mid H, s \mid v]:= \frac{\langle f\mid H, s \mid v\rangle}{(2\pi i)^{j-i-1}}.
$$

\paragraph{The map $\Phi$.} Let $\widetilde {\cal H}$ be the set of split framed mixed Hodge-Tate structures. We define first a map of sets
\[
\Phi:  \widetilde {\cal H}   \longrightarrow {\rm T}(\C).
\]
For each $0<p < n$, choose a basis $ \{v^\alpha_{p}\}$ in $\Hom(\mathbb{Q}(p), \gr^W_{-2p}H)$. Let $\{f_\alpha^{p}\}$ be the dual basis.

\bd
  \la{DEF4.1}

  We define a map $\Phi$ recursively:
 \begin{equation}
  \la{REC}
  \Phi(f^n, H,v_0; s):=-\sum_{0\leq p\leq n-1}  \sum_{\alpha}
  [ f^n \mid H,s \mid v^\alpha_p] \otimes (2 \pi i)^{\otimes(n-p-1)} \otimes \Phi(f_\alpha^p, H, v_0;s).
\end{equation}
\ed

 Definition \ref{DEF4.1}  is  independent of the choice of bases since the following element represents   the identity map on $\gr^W_{-2p}H_\mathbb{Q}$:
\be
\la{NOT}
v_p\otimes f^p := \sum_\alpha v^\alpha_p\otimes f_\alpha^p.
\ee

Below we use the following convention: whenever $f^p$ and $v_p$ appear in the same formula, this means that we take the sum over the pairs $f^p_\alpha$ and $v_p^\alpha$ of dual basis vectors, just like   in (\ref{CONV}).

 We can write the definition of $\Phi$  as a much nicer  identity:
\begin{equation}
  \la{REC2}
  \sum_{0\leq p\leq n} [ f^n \mid H,s \mid v_p]  \otimes (2 \pi i)^{\otimes(n-p-1)} \otimes \Phi(f^p, H, v_0;s)=0.
\end{equation}

Explicitly, formula (\ref{REC})  can be written schematically as
\be
     \Phi(f^n, H,v_0; s):=
    \sum_{\substack{1\leq k\leq n \\ 0=i_0<i_1<\cdots <i_k=n}}
     (-1)^k \bigotimes_{l=1}^{k} \left([ f^{i_l}\mid H,s \mid v_{i_{l-1}}] \otimes (2 \pi i)^{\otimes (i_l-i_{l-1}-1)}\right).
\ee

Given two splittings $S_W,S_W^\prime$ of $(H,v_0, f^n)$, we have
$$
S_W(v_{q})-S_W^\prime(v_{q})\in W_{-2p-2}H_\mathbb{Q}.
$$
 Therefore     $
S_{W}^{-1} \circ S_{W}^\prime$ is a lower triangular unipotent transformation with coefficients in $\mathbb{Q}$:
 $$
  S_{W}^{-1} \circ S_{W}^\prime(v_{q})=v_{q}+ \sum_{p>q} c^p v_{p}, ~~~~c^p\in \mathbb{Q}.
$$
 Let us set:
\be
\begin{split}
 &N:=S_{W}^{-1} \circ S_{W}^\prime, ~~~~[f^r\mid N \mid v_q]:=\frac{f^r \circ N \circ v_q}{(2\pi i)^{r-q+1}}\\
 \end{split}
 \ee
  The  matrix elements of the operator $N$ are rational, so following the normalization of (\ref{period}), $[f^r\mid N\mid v_q]  \in \mathbb{Q}(1)$. So we have:
  \be
\begin{split}
  \langle f^p\mid H,s^\prime\mid v_q\rangle &=f^p\circ S_{HT}\circ S_W^\prime\circ v_q\\
  &=f^p\circ S_{HT}\circ S_W\circ S_W^{-1} \circ S_W^\prime\circ v_q\\
  &=f^p\circ S_{HT}\circ S_W \circ N\circ v_q \\
  &=\sum_{ q \leq r \leq p}(f^p\circ S_{HT}\circ S_W \circ v_r)(f^r \circ N\circ v_q) \\
  &=\sum_{ q \leq r \leq p} \langle f^p\mid H,s\mid v_r\rangle \cdot ( f^r\circ N\circ v_q).
\end{split}
\ee

In other words,
\be
\la{54}
[f^p\mid H,s^\prime\mid v_q]=\sum_{ q \leq r \leq p} [f^p\mid H,s\mid v_r] \cdot \frac{[f^r\mid N\mid v_q]}{2\pi i}, ~~~~ \frac{[f^r\mid N\mid v_q]}{2\pi i}\in \mathbb{Q}.
\ee

  By the defining equation (\ref{REC2}), we have:  
\begin{align*}
  0=& \sum_{0\leq p\leq n}
   [f^n \mid H,s^\prime \mid v _p]  \otimes (2 \pi i)^{\otimes(n-p-1)} \otimes \Phi(f ^p, H, v_0;s^\prime)\\
  =& \sum_{0\leq p\leq q\leq n}
  [f^n \mid H,s \mid v _q] \cdot \frac{[f^r\mid N\mid v_q]}{2\pi i}\otimes (2 \pi i)^{\otimes(n-p-1)} \otimes
  \Phi(f^p, H, v_0;s^\prime)\\
  =& \sum_{0\leq p\leq q\leq n}
  [f^n \mid H,s \mid v _q] \otimes (2 \pi i)^{\otimes(n-q-1)} \otimes
  [f^q \mid N \mid v _p] \otimes (2 \pi i)^{\otimes(q-p-1)} \otimes
  \Phi(f^p, H, v_0;s^\prime)\\
  =& \sum_{0\leq q\leq n}
  [ f^n \mid H,s \mid v _q]  \otimes (2 \pi i)^{\otimes(n-q-1)} \bigotimes
  \sum_{0\leq p \leq q}  [f^q \mid N \mid v _p] \otimes (2 \pi i)^{\otimes(q-p-1)}  \otimes
  \Phi(f^p, H, v_0;s^\prime).
\end{align*}
Here the second and third equation follows from  (\ref{54}).

By comparing this to equation (\ref{REC2}) and by the induction, we have
\begin{equation}
  \label{Split}
  \Phi(f^n, H, v_0; s)= \sum_{0\leq p \leq n}  [f^n \mid N \mid v_p] \otimes (2 \pi i)^{\otimes(n-p-1)} \otimes
  \Phi(f^p, H, v_0;s^\prime).
\end{equation}

\bl
  Let $s, s^\prime$ be  splittings of an $n$-framed mixed Hodge-Tate structure $(f^n, H,v_0)$. Then
  $$
  \Phi(f^n,H,v_0; s)-\Phi(f^n, H,v_0; s^\prime)\in {\rm Ker}({\rm pr}).
  $$
\el

\begin{proof} We have
\[
  \Phi(f^n, H, v_0; s) -\Phi(f^n, H, v_0; s^\prime)
  =\sum_{0\leq p \leq n-1}   [f^n \mid N \mid v_p] \otimes (2 \pi i)^{\otimes(n-p-1)} \otimes
  \Phi(f^p, H, v_0;s^\prime).
\]
 Since $[f^n \mid N \mid v_p]\in \mathbb{Q}(1)$, the lemma follows.
\end{proof}

Let us set 
$$
\widetilde{\mathcal{P}}:={\rm pr}\circ \Phi.
$$ 
\bc
  \label{IndSp}
  $\widetilde{\mathcal{P}}$ does not depend on splitting.
\ec

\bl \label{IndEq}
The map $\widetilde{\mathcal{P}}$   depends only  on the equivalence class of a framed Hodge-Tate structure.
\el

\begin{proof} Let $H \hookrightarrow  H^\prime$ be an injective  map of framed mixed Hodge-Tate structures.
Choose splitting $S_W, S_W^\prime$ such that $S_W$ is the restriction of $S_W^\prime$.
Choose a basis $\{v^1_{p},\cdots, v_{p}^{ d_p^\prime}\}$ of $\Hom(\mathbb{Q}(p), \gr^W_{-2p}H^\prime)$   such that $\{v^1_{p},\cdots, v_{p}^{d_p}\}$ is a basis of $\Hom(\mathbb{Q}(p), \gr^W_{-2p}H)$.
Let $\{f^{p}_{1},\cdots, f^{p}_{d_p^\prime}\}$ be the dual basis.
Then for $j\leq d_{p}$ we have  $[ f^{p^\prime}_{j^\prime}\mid H^\prime, s^\prime \mid v_{p}^{j}] =0$   unless  $j^\prime\leq d_{p^\prime}$.
Starting with $v_0^\prime=v_0\in H$, by the induction the nonzero terms in $\Phi({f^n}^\prime, H^\prime, v_0^\prime; s^\prime)$ are exactly the terms in $\Phi(f^n, H, v_0; s)$. So $\Phi({f^n}^\prime, H^\prime, v_0^\prime; s^\prime)=\Phi(f^n, H, v_0; s)$.

By duality the claim is valid  for a surjective map $H \twoheadrightarrow H^\prime$. Since the category of mixed Hodge-Tate  structures is abelian, the general case follows.  \end{proof}

By Corollary \ref{IndSp} and Lemma \ref{IndEq}, $\widetilde{\mathcal{P}}$ descends to a well defined map $\mathcal{P}:\ \mathcal{H}_\bullet \rightarrow \mathcal{A}_\bullet$. So the map $\mathcal{P}$ is
the unique map making the following diagram commutative:
\be
\label{CommDiag}
\begin{gathered}
\xymatrix{
\widetilde {\cal H}_\bullet \ar[r]^-{\Phi} \ar[d]^{\rm p}& {\rm T}(\C) \ar[d]^{\rm pr}\\
{\cal H}_\bullet \ar[r]^-{\mathcal{P}}&{\cal A}_\bullet(\C)
}
\end{gathered}
\ee
 It is easy to see that it respects the abelian group structures, providing a map of graded $\Q$-vector spaces
 $$
 {\cal P}:\ \mathcal{H}_\bullet \rightarrow \mathcal{A}_\bullet.
 $$

 \subsection{The  map ${\cal P}$ is an isomorphism of coalgebras.}

 Although $\widetilde{\mathcal{H}}_\bullet$ does not carry a $\mathbb{Q}$-vector space structure, the definitions of the product and the coproduct are still valid if
 we  allow formal linear combinations  of elements of $\widetilde{\mathcal{H}}_\bullet$. The map $\Phi$ does not commute with the coproduct. In constrast with this, we have

   \bp \la{PROPP}The map
  ${\cal P}$ is a    homomorphism of graded coalgebras:
$$
{\cal P}: {\cal H}_\bullet \lra {\cal A}_\bullet.
$$ \ep

\begin{proof} We first show that
\be\label{CompCop}
\Delta_{\cal A} \circ {\rm pr}\circ \Phi = {\rm pr}\circ \Phi\circ \Delta_{\widetilde{\mathcal{H}}}.
\ee
Let us look at $\Delta_{\cal A} \circ {\rm pr}\circ \Phi(f^n, H, v_0; s)$. Each summand of the sum defining  $\Phi(f^n, H, v_0; s)$ is a monomial given by a
tensor product of factors, which are  either  $2\pi i$, or matrix elements $[f^q \mid H ,s \mid  v_p]$.
Applying the map ${\rm pr}$, we apply the exponent to the first factor. Then, applying the coproduct $\Delta_{\cal A}$, we apply the exponential to one of the other factors. If this factor is $2\pi i$, we get $1\in \C^*$, so the resulting contribution is zero. Therefore we need to consider only the result of application of the exponential to the matrix elements.
But this, by the very definition, coincides with the ${\rm pr}\circ \Phi\circ \Delta_{\widetilde{\mathcal{H}}}(f^n, H, v_0; s)$.

  Recall the commutative diagram (\ref{CommDiag}). Equation (\ref{CompCop}) implies:
\begin{align*}
  \Delta_\mathcal{A}\circ \mathcal{P} \circ {\rm p}
  &= \Delta_{\mathcal{A}} \circ {\rm pr} \circ \Phi\\
  &= {\rm pr}\circ \Phi\circ \Delta_{\widetilde{\mathcal{H}}}\\
  &= \mathcal{P}\circ {\rm p} \circ \Delta_{\widetilde{\mathcal{H}}}\\
  &= \mathcal{P} \circ \Delta_{\mathcal{H}} \circ {\rm p}.
\end{align*}

Since the projection ${\rm p} $ is surjective, we see that $\mathcal{P}$ commutes with $\Delta$.
\end{proof}

Proposition \ref{PROPP} implies  that the map ${\cal P}$
 gives rise to a map of the graded cobar complexes

\be \la{MCO}
\begin{gathered}
\xymatrix{
{\cal H}_\bullet\ar[r]^-{\Delta^\prime}  \ar[d]^{{\cal P}}&  \bigotimes^2 {\cal H}_\bullet\ar[d]^{\otimes^2{\cal P} }\ar[r]^{\Delta^\prime} &  \bigotimes^3 {\cal H}_\bullet\ar[d]^{\otimes^3{\cal P} } \ar[r]^-{\Delta^\prime} & \cdots \\
{\cal A}_\bullet \ar[r]^-{\Delta^\prime} &  \bigotimes^2{\cal A}_\bullet \ar[r]^{\Delta^\prime}   &\bigotimes^3{\cal A}_\bullet \ar[r]^-{\Delta^\prime} & \cdots
}
\end{gathered}
\ee


  The degree $n$ part of the cobar complex ${\rm Cobar}_n({\cal A}_\bullet)$ of the coalgebra ${\cal A}_\bullet $
 sits  in   degrees $[1,n]$. Consider the complex $\C \stackrel{\rm exp}{\lra} \C^*$   in the degrees $[0, 1]$. It is a resolution of
 $\Q(1)$.

 \bl \la{EXL}
The complex ${\rm Cobar}_n({\cal A}_\bullet)$ is isomorphic to
the complex
$$
\C^* \bigotimes \Bigl(\C \stackrel{\rm exp}{\lra} \C^*\Bigr)^{\otimes n-1}.
$$
  The complex ${\rm Cobar}_n({\cal A}_\bullet)$ is a resolution of $\C/(2\pi i)^n\Q$.\el

 \begin{proof} The  Lemma is a special case of Lemma \ref{EXL1a}. \end{proof}

For example,

\begin{align*}
  {\rm Cobar}_2({\cal A}_\bullet)
  &= \C^* \bigotimes \Bigl(\C \stackrel{\rm exp}{\lra} \C^*\Bigr).\\
{\rm Cobar}_3({\cal A}_\bullet) &= \C^* \bigotimes \Bigl(\C \stackrel{\rm exp}{\lra} \C^*\Bigr)^{\otimes 2}\\
&=\C^* \bigotimes \Bigl(\C \otimes \C \lra  \C\otimes \C^* \oplus \C^* \otimes \C \lra  \C^* \otimes \C^* \Bigr).
\end{align*}

\bt
\label{IsoCo}
The  map of complexes (\ref{MCO})  is an isomorphism:
  $$
{\cal P}^\bullet_n: {\rm Cobar}_n({\cal H}_\bullet) \lra{\rm Cobar}_n({\cal A}_\bullet).
 $$
 In particular, the map ${\cal P}_n:{\cal H}_n \lra {\cal A}_n$ is an isomorphism.
 \et

  \begin{proof}
By  (\ref{1}) and Beilinson's vanishing theorem (\ref{eq7}),  the complex  ${\rm Cobar}_n({\cal H}_\bullet) $ is a resolution of $\C/(2\pi i)^n\Q$.
 Let $\Delta_n' $ be the restriction of the restricted  coproduct $\Delta'$ to $ {\cal H}_n$.
Then we have:
$$
 {\rm Ker} (\Delta'_n)    =
\quad
 \frac{\C}{(2 \pi i)^{n}\Q}  =  {\rm Ext}_{\rm MHT_\Q}^1(\Q(0), \Q(n)).
$$

By Lemma   \ref{EXL} the same is true for the
complex ${\rm Cobar}_n({\cal A}_\bullet)$.  Moreover,  in $\mathcal{A}_\bullet$ one has:
$$
 {\rm Ker}(\Delta_n^\prime)=\mathbb{C}/2\pi i \mathbb{Q}\otimes 2\pi i \mathbb{Q} \otimes \cdots \otimes 2\pi i \mathbb{Q} \cong \mathbb{C}/(2\pi i)^n \mathbb{Q}.
$$
Therefore one easily sees that the  map  ${\cal P}^\bullet_n$  induces an isomorphism
on the only non-trivial cohomology group $\C/(2\pi i)^n\Q$. Arguing by the induction on $n$ the theorem follows.
\end{proof}

\subsection{The map ${\cal P}$ is an isomorphism of Hopf algebras}
\bp
  The map $\Phi$ is a morphism of algebras\ep

\begin{proof}
  It suffices to show $m_{s,t}\circ \Phi =\Phi\circ m_{s,t}$.
  We prove this by the induction on $s$ and $t$.

  When $s=0$ or $t=0$, the equation is trivial.

  For $s=t=1$, by the very definition of the map $\Phi$, we have:
  \begin{align*}
    \Phi\circ m (H, H^\prime)=&-[ f^1\otimes {f^\prime}^1\mid H\otimes H^\prime\mid v_0\otimes v_0^\prime] \otimes 2\pi i\\
    &+[ f^1\otimes {f^\prime}^1 \mid  H\otimes H^\prime \mid v_0\otimes v_1^\prime] \otimes [ f^0 \otimes {f^\prime}^1 \mid  H\otimes H^\prime\mid v_0\otimes v_0^\prime] \\
    &+ [ f^1\otimes {f^\prime}^1 \mid  H\otimes H^\prime \mid v_1\otimes v_0^\prime] \otimes [ f^1 \otimes {f^\prime}^0 \mid  H\otimes H^\prime\mid v_0\otimes v_0^\prime]. \\
  \end{align*}
  By elaborating the definition of the tensor product, we write this as follows:  
  \begin{align*}
    &-(2\pi i)^{-1}[f^1\mid H\mid v_0] [{f^\prime}^1\mid H^\prime\mid v_0^\prime]\otimes 2\pi i\\
    &+(2\pi i)^{-1}[f^1\mid H\mid v_0][{f^\prime}^1 \mid  H^\prime \mid  v_1^\prime]
    \otimes (2\pi i)^{-1}[ f^0\mid H\mid v_0][{f^\prime}^1 \mid  H^\prime\mid  v_0^\prime] \\
    &+(2\pi i)^{-1}[f^1\mid H\mid v_1][{f^\prime}^1 \mid  H^\prime \mid  v_0^\prime]
    \otimes (2\pi i)^{-1}[ f^1\mid H\mid v_1][{f^\prime}^0 \mid  H^\prime\mid  v_0^\prime]. \\
  \end{align*}
  We observe that four of the factors above are equal to $2 \pi i$:
  $$
  [ {f^\prime}^1 \mid  H^\prime \mid  v_1^\prime]= [ f^0\mid H\mid v_0] =  [ f^1\mid H\mid v_1]= [ {f^\prime}^0 \mid H^\prime\mid  v_0^\prime] = 2 \pi i.
  $$   Therefore we get:
  \begin{align*}
    &-(2\pi i)^{-1}[f^1\mid H\mid v_0] [{f^\prime}^1\mid H^\prime\mid v_0^\prime] \otimes 2 \pi i\\
    &+[ f^1\mid H\mid v_0]\otimes  [ {f^\prime}^1 \mid  H^\prime\mid  v_0^\prime] \\
    &+[{f^\prime}^1 \mid  H^\prime \mid  v_0^\prime] \otimes   [ f^1\mid H\mid v_1] \\
    =&m (\Phi (H), \Phi(H^\prime)).
  \end{align*}
  Assume the equation is proved for $(a,b)$ such that $a\leq s, b\leq t, a+b \leq s+t-1$. Then
  \[\Phi \circ m_{s,t} (H, H^\prime) = \Phi (f^s\otimes {f^\prime}^t, H\otimes H^\prime, v_0\otimes v_0^\prime).\]
  Using the inductive definition (\ref{REC}) of the product, we write this as follows:
  \[-\sum_{\substack{0\leq p+p^\prime\leq s+t-1   }}
    [ f^s\otimes {f^\prime}^t \mid H\otimes H^\prime \mid v_p \otimes v^\prime_{p^\prime}] \otimes (2\pi i)^{\otimes (s+t-p-p^\prime-1)}
    \bigotimes \Phi (f^p\otimes {f^\prime}^{p^\prime}, H\otimes H^\prime, v_0\otimes v_0^\prime).
  \]

  Using the definition of the tensor product in the left factor,  and the inductive assumption in the right factor, we get:
  \be
  \bs
 & -\sum_{\substack{0\leq p+p^\prime\leq s+t-1  }}
   (2\pi i)^{-1} [f^s\mid H \mid v_p] [{f^\prime}^t \mid H^\prime \mid v^\prime_{p^\prime}] \otimes (2\pi i)^{\otimes (s+t-p-p^\prime-1)}
    \bigotimes \\
    &m(\Phi (H, v_0, f^p),\Phi({f^\prime}^{p^\prime}, H^\prime, v_0^\prime)).\\
\end{split}
\ee
  Let us split this sum into three terms: in the first term we have  $p=s$, in the second $p'=t$, and the last one is the rest, $p<s, p'<t$:
  \begin{align*}
    &-\sum_{\substack{0\leq p \leq s-1 }}     [ f^s\mid H \mid v_p] \otimes (2\pi i)^{\otimes (s-p-1)}
    \bigotimes m(\Phi (H, v_0, f^p),\Phi({f^\prime}^t, H^\prime, v_0^\prime))\\
    &-\sum_{\substack{0\leq p^\prime\leq t-1  }}     [ {f^\prime}^t \mid H^\prime \mid v^\prime_{p^\prime}] \otimes (2\pi i)^{\otimes (t-p^\prime-1)}
    \bigotimes m(\Phi (H, v_0, f^s),\Phi({f^\prime}^{p^\prime}, H^\prime, v_0^\prime))\\
    &-\sum_{\substack{0\leq p \leq s-1,\ 0\leq p^\prime\leq t-1  }}
  (2\pi i)^{-1}[ f^s\mid H \mid v_p] [ {f^\prime}^t \mid H^\prime \mid v^\prime_{p^\prime}]   \otimes (2\pi i)^{\otimes (s+t-p-p^\prime-1)}
    \bigotimes\\
    &m(\Phi (f^p, H, v_0),\Phi({f^\prime}^{p^\prime}, H^\prime, v_0^\prime)).
  \end{align*}
  Using the inductive definition (\ref{REC}) of the product $m$, we recognize in the sum of the three terms above the definition of the product of the two terms below, where we use the sign $*$  for the product   to make the notation compact:
  \begin{align*}
    &\Bigl(-\sum_{\substack{0\leq p \leq s-1    }}  [ f^s\mid H \mid v_p] \otimes (2\pi i)^{\otimes (s-p-1)}\otimes \Phi (f^p, H, v_0)\Bigr) * \\
    &\Bigl( -\sum_{\substack{0\leq p^\prime\leq t-1   }}  [ {f^\prime}^t \mid H^\prime \mid v^\prime_{p^\prime}]  \otimes (2\pi i)^{\otimes (t-p^\prime-1)}\otimes \Phi({f^\prime}^{p^\prime}, H^\prime, v_0^\prime)\Bigr)\\
    =&m_{s,t}\circ \Phi(H, H^\prime).
  \end{align*}

  So $m_{s,t}\circ \Phi =\Phi\circ m_{s,t}$.
  Therefore $\Phi$ is a homomorphism of algebras.
\end{proof}

\bc
\label{IsoAl}
The map $\mathcal{P}$ is a morphism of algebras.
\ec

By Theorem \ref{IsoCo} and Corollary \ref{IsoAl}, we have
\bt \la{MAPPI}
The map $\mathcal{P}:\ \mathcal{H}_\bullet \rightarrow \mathcal{A}_\bullet$ is an isomorphism of graded Hopf algebras.
\et

\paragraph{The inverse map.}

We define   a map $\eta:\ \mathcal{A}_\bullet \longrightarrow \mathcal{H}_\bullet$  by setting

\be \la{ETA}
\begin{split}
&\eta: A_1\otimes a_2 \otimes\cdots \otimes a_n \in \C^* \otimes \C^{\otimes n-1}\lms \\
&  \begin{pmatrix}
     1\\
     & 2\pi i\\
     & & & \ddots\\
     & & & & (2\pi i)^{n-1}\\
     & & & & & (2\pi i)^n
   \end{pmatrix} \cdot
   \begin{pmatrix}
     1\\
     a_{n} & 1&\\
     & \ddots & \ddots & \\
     & & a_2 & 1&\\
     & & &  \log A_1 & 1
   \end{pmatrix}\\
 \end{split}
\ee

 For example:

 $\eta(A_1 \otimes a_2 )=\begin{pmatrix}
  1\\
      & 2\pi i &\\
     &     & (2\pi i)^2
     \end{pmatrix} \begin{pmatrix}
  1\\
     a_{2} & 1&\\
     &   \log A_{1}  & 1
     \end{pmatrix} = \begin{pmatrix}
  1\\
  a_{2} & 2\pi i&\\
    & 2\pi i \log A_{1} & (2\pi i)^2&
     \end{pmatrix}$

     \bp The map $\eta$ the inverse to the map ${\cal P}$.
     \ep

     \begin{proof} The map $\eta$ does not depend on the choice of $\log (A_1)$. It is additive with respect to $\log A_1, a_2, ..., a_n$. So it is a well defined map of sets.
     It is evident that ${\cal P}\circ \eta = {\rm Id}$: all periods of length $>1$ of the mixed Hodge-Tate structure with the period matrix (\ref{ETA}) are zero.
     Since  the map ${\cal P}$ is an isomorphism by Theorem \ref{MAPPI}, the claim follows.
     \end{proof}

\subsection{An application: the big period maps} \label{App}

 The big period map   is    the following group homomorphism,   constructed in \cite{G96}, \cite{G15}:
 \[P_n:\ \mathcal{H}_n \rightarrow \mathbb{C}^* \otimes_\mathbb{Q} \mathbb{C}(n-2).\] The map $P_2$ for the dilogarithm ${\rm Li}_2(x)$
 goes back to S. Bloch \cite{Bl77}, \cite{Bl78}.

Let  us define  a similar   homomorphism of abelian groups for the Hopf algebra $\mathcal{A}_\bullet$:
\be
\begin{split}
&P_{n,\mathcal{A}}:\ \mathcal{A}_n \rightarrow \mathbb{C}^*\otimes_\mathbb{Q} \mathbb{C}(n-2).\\
&P_{n,\mathcal{A}}(A_1\otimes a_2\otimes\cdots \otimes a_n)=-A_1\otimes  \frac{a_2 \cdot \ldots  \cdot a_k}{(2\pi i)^{n-2}} \otimes (2\pi i)^{n-2}.\\
\end{split}
\ee
  Then, by the very definition, the map  $P_n$ factors through $\mathcal{A}_\bullet$, i.e. there is a commutative diagram with $P_n=P_{n, \mathcal{A}}\circ {\cal P}$:

  \centerline{
    \xymatrix{
      \mathcal{H}_n \ar[rr]^{\mathcal{P}} \ar[rd]_{P_n}& & \mathcal{A}_n \ar[ld]^{P_{n, \mathcal{A}}}\\
       & \mathbb{C}^* \otimes_{\mathbb{Q}} \mathbb{C}(n-2)&
    }
  }

In \cite{G15},   another version of the big period map was constructed:
\[P_n^\prime:\ \mathcal{H}_n \rightarrow \mathbb{C} \otimes_\mathbb{Q} \mathbb{C}.\]
The map  $P_n$ is the composition of $P_n^\prime$ with the map
\be
\begin{split}
&\mathbb{C} \otimes_\mathbb{Q} \mathbb{C} \rightarrow \mathbb{C}^* \otimes_\mathbb{Q} \mathbb{C}(n-2)\\
&a \otimes b \mapsto \exp(2\pi i a)\otimes 2\pi i b \otimes (2\pi i)^{n-2}.\\
\end{split}
\ee

We  define
\be
\begin{split}
&P_{n,\mathcal{A}}^\prime:\ \mathcal{A}_n \rightarrow \mathbb{C}\otimes_\mathbb{Q} \mathbb{C}\\
&A_1\otimes a_2\otimes\cdots \otimes a_n \mapsto
-\frac{\log A_1 }{2\pi i} \otimes \frac{a_2\cdot \ldots  \cdot a_k}{(2\pi i)^{n-1}} + 1\otimes \frac{\log A_1 \cdot a_2 \cdot \ldots \cdot a_k}{(2\pi i)^{n}}.\\
\end{split}
\ee
Then there is a commutative diagram with $P_n^\prime=P_{n,\mathcal{A}}^\prime \circ {\cal P}$:

\centerline{
  \xymatrix{
    \mathcal{H}_n \ar[rr]^{\mathcal{P}} \ar[rd]_{P_n^\prime}& & \mathcal{A}_n\ar[ld]^{P_{n,\mathcal{A}}^\prime}\\
    &\mathbb{C} \otimes_\mathbb{Q} \mathbb{C}&
  }
}

\section{The Galois group of the category of variations of Hodge-Tate structures } \la{SEC4}

The graded Hopf algebra $H^0_{\cal D}({\cal A}_\bullet(\Omega))$ looks as follows:
\be
\bs
&H^0_{\cal D}({\cal A}_n(\Omega)) = \{x \in {\cal A}_n({\cal O})~|~ {\cal D}(x)=0 \} \subset {\cal O}^* \otimes {\cal O}^{\otimes n-1}.\\
\end{split}
\ee

For example, one has
\be
\bs
&H^0_{\cal D}({\cal A}_1(\Omega)) = {\cal O}^*.\\
&\\
&H^0_{\cal D}({\cal A}_2(\Omega)) = {\rm Ker} \Bigl({\cal O}^* \otimes {\cal O}\stackrel{{\cal D}}{\lra} \Omega^1\Bigr),\\
&{\cal D}(F \otimes g)= d\log F \cdot g.\\
&\\
&H^0_{\cal D}({\cal A}_3(\Omega )) = {\rm Ker} \Bigl({\cal O}^* \otimes {\cal O}\otimes {\cal O}\stackrel{{\cal D}}{\lra} \Omega^1\otimes {\cal O} \bigoplus  {\cal O}^* \otimes \Omega^1  \Bigr), \\
&{\cal D}(F \otimes g_1 \otimes g_2)= d\log F \cdot g_1 \otimes g_2 + F \otimes d\log g_1 \cdot g_2.\\
\end{split}
\ee

Denote by ${\cal H}_\bullet(X)$ the sheaf of Hopf algebras of framed variations of $\Q$-Hodge-Tate structures on   $X$.
Restricting a framed variation to a point $x \in X$ we get a map of Hopf algebras
$$
r_x: {\cal H}_\bullet(X) \lra {\cal H}_\bullet.
$$
Applying pointwise the map ${\cal P}$, that is considering the composition ${\cal P}\circ r_x$, we get an injective  map
\be \la{POIN}
{\cal P}_X: {\cal H}_\bullet(X) \hra {\cal A}_\bullet({\cal O}).
\ee

\bl The Griffith transversality  implies that the  image of the map (\ref{POIN}) is killed by the differential ${\cal D}$, providing us with a map
\be \la{MAPPP}
{\cal P}_X: {\cal H}_\bullet(X) \hra \{x \in {\cal A}_\bullet({\cal O})~|~ {\cal D}(x)=0 \} \subset {\cal A}_\bullet({\cal O}).
\ee
\el

\begin{proof}
  Let
  \begin{align*}
  &D_1: X_1 \otimes x_2 \otimes\cdots \otimes x_n \mapsto d\log X_1 \cdot x_2 \otimes x_3 \otimes \cdots \otimes x_n;\\
  &D_i: X_1 \otimes x_2 \otimes \cdots \otimes x_n \mapsto X_1 \otimes x_2 \otimes\cdots \otimes d x_i \cdot x_{i+1}\otimes \cdots \otimes x_n, \qquad i>1.
  \end{align*}
  Then $D$ is a signed sum of $D_i$. It suffices to show $D_i(\mathcal{P}_X (H))=0$.

  The terms   in $\mathcal{P}_X (H)$ where the $i$-th term is   a constant are killed by $D_i$. So it remains to consider  the terms with $f^{n-i+1}$ only. At the $i$-th and $(i+1)$-th tensor factor, it is either $[f^{n-i+1} \mid H \mid v_{n-i}]\otimes [f^{n-i}\mid H \mid v_p]$ or $[f^{n-i+1} \mid H \mid v_p]\otimes 2\pi i$, where $p \leq n-i-1$. By the parts iii),  iv) of  \cite[Theorem 2.8]{G15},
  the  Griffith transversality implies
  \be
  \label{Gt}
 2\pi i  \cdot d [f^{n-i+1} \mid H \mid v_p]  = d [f^{n-i+1} \mid H \mid v_{n-i}]\cdot  [f^{n-i}\mid H \mid v_p].
  \ee
  This means that the map $D_i$ cancels them  out.
\end{proof}

\bt  \la{MAPPP2} The map (\ref{MAPPP}) is an isomorphism of sheaves of graded Hopf algebras on $X$:
\be \la{MAPPP1}
{\cal P}_X: {\cal H}_\bullet(X) \stackrel{\sim}{\lra} H^0_{\cal D}({\cal A}_\bullet(\Omega )) = \{x \in {\cal A}_\bullet({\cal O})~|~ {\cal D}(x)=0 \}.
\ee
 \et

\begin{proof}

Let $t\in \mathcal{A}_n(\mathcal{O})$ be an element satisfying  ${\cal D}(t)=0$. Let us
   show that there exists a variation of framed mixed $\Q$-Hodge-Tate structure $M$ whose image under
  the map $\mathcal{P}_X$ is $t$.

  We prove this by induction.
  The $n=1$ case is trivial.

  For $n=2$, let $t=\sum X_i \otimes y_i$. Then $t$ defines  pointwise a collection $\{M_x\}$ of framed mixed Hodge-Tate structures  on $X$. We may assume that
${\rm gr}^W_ 0 M_x= \Q(0)$, ${\rm gr}^W_ 4M_x=\Q(2)$, and ${\rm gr}^W_ {2k} M_x= 0$ unless $k=0, 1,2$.   Then the Griffith transversality condition boils down by (\ref{Gt})
  \[ 2\pi i \cdot d [f^2 \mid H \mid v_0] = d [f^2 \mid H \mid v_1]\cdot  [f^1\mid H \mid v_0].\]
  But this condition just means that  ${\cal D}(t)=0$.

  For $n\geq 3$, assume that the map $\mathcal{P}_{X, < n}:\ \mathcal{H}_{<n}(X) \rightarrow H^0_{\mathcal{D}}(\mathcal{A}_{<n} (\Omega))$ is an isomorphism.
  Let us pick  an element $t\in \mathcal{A}_n(\mathcal{O})$ such that $ \mathcal{D}(t)=0$.

  Step 1: Let $t= \sum_k X_k \otimes x_k$, where $X_k \in \mathcal{A}_{n-1}(\mathcal{O})$, $x_k\in \mathcal{O}$.
  We may assume $x_k$'s are linearly independent. Indeed,  if $x_1= \sum_{k\geq 2} c^k x_k$, then we rewrite $t$ as $\sum_{k\geq 2} (X_k + c^k X_1) \otimes x_k$.

  Since $x_k$'s are linearly independent, ${\cal D}(X_k \otimes x_k)=0$ implies ${\cal D}(X_k)=0$.

  Step 2: Since ${\cal D}(X_k)=0$, by the inductive assumption  there exists a split variation of framed mixed $\Q$-Hodge-Tate structures $M_k$ such that $\mathcal{P}_X(M_k)=X_k$.

  We need the following technical lemma.

\bl
  \label{LI}
  Let $\{M_k\}$ be a finite set of variations of $n$-framed mixed Hodge-Tate structures. Then there exists a variation of mixed Hodge-Tate structure $M$ with a basis $\{v^a_i\}$ of $\gr^W_{-2i}M$, $i=0, ..., n$, such that $(f^n, M, v_0^k)\sim M_k$, and $\mathcal{P}_X( f^n, M , v^a_i)$ are linearly independent for any $1\leq i\leq n-1$.
\el

\begin{proof}
  Let $S$ be the set of splitted variations of mixed Hodge-Tate structures with a basis $\{v^a_i\}$ of $\gr^W_{-2i}$ such that $(f^n, M, v_0^k)\sim M_k$. It is non-empty since $\bigoplus M_k$ with $f^n:=\sum f^n_k$ and $v_0^k$ is in $S$. Let $M$ be an element in $S$ with the lowest dimension. We will show by induction that $\mathcal{P}_X( f^n, M , v^a_i)$'s are non-zero and linearly independent for any $1\leq i\leq n-1$.

  We need the following lemma to reduce the dimension of $M$:

  \bl
    \label{Reduce}
    Let $M$ be a variation of mixed Hodge-Tate structures with a splitting and a basis $\{v_i\}$ of $\gr^W_{-2i} M$. If for some $v_i^a$, $[f^j \mid M \mid v_i^a]=0$ for any $j>i$ and $f^j$, then the image of $S_W \circ v_i^a$, still denoted by $v_i^a$, is a subvariation of mixed Hodge-Tate structures.

    Let $M^\prime=M/v_i^a$. Then $M^\prime$ has the same periods as $M$, i.e., for any $(j,b), (k, c)\neq (i,a)$, $[f^j_b \mid M^\prime \mid v_k^c]=[f^j_b \mid M \mid v_k^c]$.
  \el

  \begin{proof}
    $\langle f^j \mid M \mid v_i^a\rangle = f^j \circ S_{HT} \circ S_W \circ v_i^a$ is the following map:
    \[\mathbb{Q}(i) \otimes \mathbb{C} \rightarrow \gr ^W_{-2i} H_\mathbb{Q} \otimes \mathbb{C} \rightarrow W_{-2i} H_\mathbb{Q} \otimes \mathbb{C} \rightarrow \bigoplus_{p\geq i} F^{-p} H_\mathbb{C}\cap W_{-2p} H_\mathbb{C} \rightarrow \bigoplus_{p\geq i} \gr^W_{-2p} H_\mathbb{C}.\]

    By the given condition, the component for $p>i$ is zero. Then the image in $H_\mathbb{Q}$ and its complexification, denoted by $v_i^a$, satisfies $W_{-2i}v_i^a=v_i^a$, $W_{-2i-1}v_i^a=0$, $F^{-i}v_i^a=v_i^a$, $F^{-i+1}v_i^a=0$. Hence $v_i^a$ is a subvariation of mixed Hodge-Tate structures. The second half of the lemma is trivial.
  \end{proof}

  Let $i=n-1$. If $\mathcal{P}_X( f^n, M , v^a_{n-1})=0$ for some $v^a_{n-1}$, then by Lemma \ref{Reduce} $M/v^a_{n-1} \in S$ has lower dimension. If $\mathcal{P}_X( f^n, M , v_{n-1}^0)= \sum_ac_a\mathcal{P}_X( f^n, M , v_{n-1}^a)$, then replace $v_{n-1}^0$ by $v_{n-1}^0- \sum c_a v_{n-1}^a$ in the basis. It is an element in $S$ with $\mathcal{P}_X( f^n, M , v_{n-1}^0- \sum c_a v_{n-1}^a)=0$, so does not have the lowest dimension.

  Let us assume that  the claim is proved for $i+1, \cdots , n-1$. Suppose that $\mathcal{P}_X( f^n, M , v^a_i)=0$ for some $v^a_i$.

  For any $i <j\leq n$, we may assume $[ f^j_{b} \mid M \mid v^a_i] \not \in \mathbb{Q}(1)-\{0\}$ for any $f^j_b$ in the basis. Otherwise we change the splitting to make it zero.
  Moreover, we may assume for any $j$, $2\pi i$ is not a rational linear combination of $\{[f^j_b\mid M\mid v^a_i]\}$.

  Then we have
  \begin{align*}
    0 &= \mathcal{P}_X( f^n, M , v^a_i)\\
    &= \sum_{j=i+1}^{n-1} \mathcal{P}_X(f^n, M, v_j) \otimes [f^j \mid M \mid v^a_i] \otimes (2\pi i)^{j-i+1}\\
    &\ \ + \exp([f^n \mid M \mid v^a_i]) \otimes (2\pi i)^{n-i+1}.
  \end{align*}

  Notice that every summand in the equation ends with a tensor factor $(2 \pi i)$, except the terms in $\mathcal{P}_X(f^n, M, v_{n-1}) \otimes [f^{n-1} \mid M \mid v^a_i]$. With the assumption that $2\pi i$ is not a linear combination of $[f^{n-1} \mid M \mid v^a_i]$, we see $\mathcal{P}_X(f^n, M, v_{n-1}) \otimes [f^{n-1} \mid M \mid v^a_i]=0$. For the same reason, we can show inductively that $\mathcal{P}_X(f^n, M, v_j) \otimes [f^j \mid M \mid v^a_i] \otimes (2\pi i)^{j-i+1}=0$ for all $i<j<n$ and $\exp([f^n \mid M \mid v^a_i]) \otimes (2\pi i)^{n-i+1}=0$.

  Since for any $i < j < n$, $\mathcal{P}_X(f^n, M, v_j)$'s are linearly independent by inductive assumption, the equations above shows $[f^j \mid M \mid v^a_i]=0$ for any $f^j$. We also assumed $[f^n \mid M \mid v^a_i]\not \in \mathbb{Q}(1)-\{0\}$, so $\exp([f^n \mid M \mid v^a_i])=0$ implies $[f^n \mid M \mid v^a_i]=0$. Therefore $[f^j \mid M \mid v^a_i]=0$ for any $i<j\leq n$ and any $f^j$ in the chosen basis.

  By Lemma \ref{Reduce}, $M/v^a_i \in S$ has a smaller dimension. Contradiction. Therefore $\mathcal{P}_X( f^n, M , v^a_i)$ is non-zero.

  By similar argument as in $i=n-1$, $\{\mathcal{P}_X(f^n, M , v^a_i)\}$ are linearly independent.
\end{proof}

Let $M^\prime$ with $\{v_i^{a\prime}\}$ satisfies Lemma \ref{LI} for $\{M_k\}$. Then let us construct an $M\in \mathcal{H}_n(X)$ as follows.
Consider a $\Q$-vector space  $M=\mathbb{Q}(0)\oplus M^\prime(1)$. Choose a splitting and a basis $\{v^a_i\}$ on $M^\prime(1)$ provided by   $M^\prime$.   We extend it to
 a splitted    mixed Hodge-Tate structure on $M$ such that
 $[f^1_k \mid M \mid v_0]=x_k$, and $[f^i \mid M \mid v_0]=0$ for $i\geq 2$.

  Obviously $\mathcal{P}_X(M)= t$. Then ${\cal D}_{n-1}(t)=0$ implies
  \[\mathcal{P}_X(f^n, M, v_2) \otimes d [f^2 \mid M \mid v_1]\cdot  [f^1\mid M \mid v_0] =0.\]
  Since $\mathcal{P}_X(f^n, M, v_2)=\mathcal{P}_X(f^n(-1), M^\prime, v_2(-1))$, by Lemma \ref{LI} they are linearly independent, hence
  \[d [f^2 \mid M \mid v_1]\cdot  [f^1\mid M \mid v_0]=0.\]
This implies (\ref{Gt}).
  Therefore the Griffith transversality condition is satisfied, and $M$ is a variation of mixed $\Q$-Hodge-Tate structures, with $\mathcal{P}_X(M)=t$. This shows $\mathcal{P}_X$ is surjective.
  \end{proof}

 \bc \la{CCORR}
 Given a complex disc $U$, the category of graded comodules over the commutative Hopf algebra $H^0_{\cal D}({\cal A}_\bullet(\Omega_U ))$ is
 canonically equivalent to the category of variations of mixed Hodge-Tate structures on $U$.  \ec
 
 \begin{proof} The category ${\rm VMHT}_{U}$ of variations of mixed Hodge-Tate structures on a complex manifold $U$ is a mixed Tate $\Q$-category,  
 with the Tate object given by the constant variation $\Q(1)_U$. Thanks to the basics  about mixed Tate categories, see \cite{BD92} or Appendix A, equivalence (\ref{eqcc}), 
 the category ${\rm VMHT}_{U}$ 
 is equivalent to the category of graded comodules  over the corresponding Hopf algebra ${\cal H}_\bullet(U)$. So 
 Corollary \ref{CCORR}  follows from this and Theorem \ref{MAPPP2}. 
 \end{proof}

   Therefore we arrive at a variant of Theorem \ref{TH5.1}.
\bc i) The weight $n$ part of the cobar complex of sheaves ${\rm Cobar}_n\Bigl(H^0_{\cal D}({\cal A}_\bullet(\Omega )) \Bigr)$  is quasi-isomorphic to the weight $n$ Deligne complex of sheaves 
of $X$.

ii) The  cobar complex of the Hopf algebra $H^0_{\cal D}({\cal A}_\bullet(\Omega )) $  is commutative dg-algebra in the category of  sheaves, providing another dg-module for the Deligne cohomology of $X$.
 \ec

\section{Period morphism for p-adic   mixed Tate Galois representations} \la{Sec5}

\subsection{Hopf algebras arising from Tate $\varphi$-algebras} \la{sec5.1}

We start with a variant of the notion of  Tate algebra, which we call {\it   Tate $\varphi$-algebra}. 
A slight modification of our main construction  assigns to   any   Tate $\varphi$-algebra $R$ a graded Hopf algebra ${\cal A}^\varphi_\bullet(R)$. 
Our main example of a Tate $\varphi$-algebra is given by a subalgebra  of Fontaine's ring $\B_{\rm crys}$.

\bd A    Tate $\varphi$-algebra $R$ over a field  $k$ is  a graded $k$-algebra     with an invertible  Tate line  $k(1)$ in the degree $1$:
\be \la{invtl}
R = \bigoplus_{n \in \Z} R_n, ~~~~k(1)-\{0\}  \subset R^*_1.
\ee

 \ed

For any  $n\in \Z$ set  $k(n) := k(1)^{\otimes n}$. Then there is      an isomorphism
 $R_1  \otimes k(n-1) \stackrel{\sim}{\lra} R_n$.  
 Set
$$
\overline {R_n}:= R_n/k(n).
$$
Consider a graded $k$-vector space
\be \la{phial}
\begin{split}
  &{\cal A}^\varphi_\bullet(R):=   \bigoplus_{n=0}^\infty {\cal A}^\varphi_n(R), ~~~~{\cal A}^\varphi_0(R) = k.\\
  &{\cal A}^\varphi_d(R):=        \overline {R_1} \otimes_{k} \underbrace{{R_1} \otimes_{k}
\ldots \otimes_{k} {R_1}}_{\text{ $n-1$ factors }}, ~~n>0.\\
\end{split}
\ee

\bt \la{Th2.2v} Let $R$ be a   Tate $\varphi$-algebra. Then the 
 graded $k$-vector space  ${\cal A}^\varphi_\bullet(R)$ has a graded  Hopf algebra structure. It is commutative if and only if $R$ is commutative.
\et

\begin{proof}
The coproduct is given by the deconcatenation, followed by the projection of the first factor to $\overline R_1$. The product
is given by the inductive formula (\ref{EQPin}) with $x_i, y_i \in R_1$. Observe that if $x_1, y_1 \in R_1$, then $x_1y_1/t \in R_1$, as
well as $t\in R_1$, so   formula (\ref{EQPin}) does make sense. The rest is identical to the proof of Theorem \ref{TH2.2}. \end{proof} 

We assume below that $R$ is commutative.
Then there is a tensor category ${\cal M}(R)$ of graded ${\cal A}^\varphi_\bullet(R)$-comodules. Denote by $k(n)$ the one dimensional comodule in the degree $n$. 
Then
\be
\begin{split}
&{\rm Ext}^1_{{\cal M}(R)}(k(0), k(n)) =  \overline R_n . \\
&{\rm Ext}^i_{{\cal M}(R)}(k(0), k(n)) =0, ~~i>1.\\
\end{split}
\ee
Indeed, in analogy with Lemma \ref{EXL1a},
the complex ${\rm Cobar}_n({\cal A}^\varphi_\bullet(R))$ is isomorphic to
the complex
$$
\overline R_1\bigotimes \Bigl(R_1 \stackrel{}{\lra} \overline R_1 \Bigr)^{\otimes n-1}.
$$
  Therefore the complex ${\rm Cobar}_n({\cal A}^\varphi_\bullet(R))$ is a resolution of $\overline R_1 \otimes k(n-1) = \overline R_n $.

 .

\subsection{The period morphism ${\cal P}_{\rm crys}$ for crystalline mixed Tate Galois modules} \la{sec5.3}

\paragraph{Tate $\varphi$-algebra ${\rm R}_{\rm crys}$ and related    graded Hopf algebra.}  As   explained in Section \ref{SECC1.1a}, there is a subalgebra ${\rm R}_{\rm crys}\subset {\rm B}_{\rm crys}$ which has 
a structure of a  Tate $\varphi$-algebra over $\Q_p$:
$$
{\rm R}_{\rm crys}= \bigoplus_{n\in \Z} {\rm B}_{\rm crys}^{\varphi =p^n}, ~~~~\Q_p(1) := \Q_p t \subset {\rm B}_{\rm crys}^{\varphi =p}.
$$   
Recall that   for each $n \in \Z$ we have the fundamental exact sequence:
$$
0 \lra \Q_pt^n \lra {\rm B}_{\rm crys}^{\varphi =p^n} \lra {\rm B}_{\rm dR}/t^n{\rm B}^+_{\rm dR} \lra 0.
$$
So in this case we have
$$
 {R_n} = {\rm B}_{\rm crys}^{\varphi =p^n}, ~~~~\overline {R_n} = {\rm B}_{\rm dR}/t^n{\rm B}^+_{\rm dR}.
$$
Applying our construction, we get    the following $\Q_p$-vector spaces
\be
\begin{split}
&{\cal A}^\varphi_n({\rm R}_{\rm crys}):=    
  {\rm B}_{\rm dR}/t{\rm B}^+_{\rm dR} \otimes_{\Q_p} \underbrace{{\rm B}_{\rm crys}^{\varphi =p}
 \otimes _{\Q_p}   \ldots \otimes _{\Q_p}  {\rm B}_{\rm crys}^{\varphi =p}}_{\text{ $n-1$ factors }}.\\
\end{split}
\ee
So we arrive at a positively graded commutative Hopf algebra over $\Q_p$:
$$
{\cal A}^\varphi_\bullet({\rm R}_{\rm crys}):= \Q_p\oplus \bigoplus_{n=1}^\infty {\cal A}^\varphi_n({\rm R}_{\rm crys}).
$$
Therefore we get   the category $ {\cal M}({\rm B}_{\rm crys})$ of graded ${\cal A}^\varphi_\bullet({\rm R}_{\rm crys})$-modules with
\be \la{CEG}
{\rm Ext}^1_{{\cal M}({\rm R}_{\rm crys})}(\Q_p(0), \Q_p(n)) = {\rm B}_{\rm dR}/t^n{\rm B}^+_{\rm dR} =  {\rm B}_{\rm crys}^{\varphi =p^n} / \Q_pt^n, ~~~~n>0.
\ee

\paragraph{Mixed Tate crystalline representations.}  Let $K$ be a finite extension of $\Q_p$. Let $\G_{K}= {\rm Gal}(\overline K/K)$ be the Galois group of $K$. 
We start with the notion of   mixed Tate Galois modules. 

\bd A p-adic $\G_K$-representation $V$ is {\it mixed Tate} if $V$ has   an increasing filtration $W_\bullet$, called the weight filtration, 
such that   ${\rm gr}^W_{-2d}V = \oplus {\Q}_p(d)$ and  
${\rm gr}^W_{-2d+1}V=0$.
\ed

 Let $K_0\subset K$ be the maximal   subfield of 
$K$ unramified over $\Q_p$. Let $V$ be a   $\G_{K}$-representation   over $\Q_p$. 
 Consider a $K_0$-vector space 
$$
{\rm D}_{\rm crys}(V) := \Bigl(V \otimes_{\Q_p} {\rm B}_{\rm crys}\Bigr)^{\G_{K}}.
$$
Recall  that  $V$ is {\it crystalline} if  ${\rm dim}{\rm D}_{\rm crys}(V) = {\rm dim}(V)$, i.e. the   map 
$
i: {\rm D}_{\rm crys}(V) \hra V \otimes_{\Q_p} {\rm B}_{\rm crys} 
$ 
  induces an isomorphism  
\be \la{CRYSTI}
i_{{\rm B}_{\rm crys}}: {\rm D}_{\rm crys}(V) \otimes_{K_0}  {\rm B}_{\rm crys}  \stackrel{\sim}{\lra}  V\otimes_{\Q_p} {\rm B}_{\rm crys} .
\ee

The  $K_0$-vector space ${\rm D}_{\rm crys}(V)$  is equipped with a $\sigma$-linear  operator $\varphi$, where $\sigma \in {\rm Gal}(K_0/\Q_p)$ is the Frobenious map, 
i.e. it is a $\varphi$-module over $K_0$. 
For example, 
\be \la{91}
{\rm D}_{\rm crys}(\Q_p(i)) = K_0[-i]
\ee 
is a one-dimensional $K_0$-vector space   with the map $\varphi$ acting by $p^{-i}\sigma$,

Given a $\varphi$-module $M$ over $K_0$, consider the {\it slope filtration} on $M$ defined as follows. Pick a lattice $M_0 \subset M$, so that  $M_0 \otimes_{{\cal O}_{K_0}} K_0 =  M$. For each integer $k$ consider 
$$
{\cal F}^kM:= \{x \in M ~|~ (p^{k}\varphi)^n (x) \mbox{ is $M_0$-bounded as $n \to \infty$}\}.
$$
It is a decreasing filtration on $M$. A different choice of $M_0$ leads to the same filtration. 

In particular, we get a {\it slope filtration} ${\cal F}^k{\rm D}_{\rm crys}(V)$. 
Below we use it, just like we use the Hodge filtration in the Hodge-Tate case, to define the p-adic period operator.

Let us denote by ${\cal M}_{\rm \G_K, crys}$ the category of p-adic mixed Tate crystalline $\G_K$-representations. It is a mixed Tate category over $\Q_p$. So it is described as the category of finite dimensional  graded comodules its fundamental graded Hopf algebra, denoted by   
${\cal H}_\bullet({\cal M}_{\rm \G_K, crys})$, see Appendix A.

   \paragraph{The p-adic crystalline period operator.}

Let $V$ be a mixed  Tate crystalline $\G_K$-representation over $\Q_p$.   Then there is a decreasing slope filtration on 
${\cal F}^\bullet$ on $ {\rm D}_{\rm crys}(V)$, and an increasing weight filtration $W_\bullet$ on $V$. Since  $V$ is crystalline, there is an isomorphism 
(\ref{CRYSTI}). 
Thanks to (\ref{91}),   the slope filtration 
${\cal F}^p {\rm D}_{\rm crys}(V) \otimes_{\Q_p} \B_{\rm crys}$ is   opposite to the weight filtration $W_\bullet V \otimes_{K_0} \B_{\rm crys}$. 
Therefore these two filtrations split $V \otimes_{\Q_p} \B_{\rm crys}$:
\be \la{94}
\begin{split}
&V \otimes_{\Q_p} \B_{\rm crys} =\bigoplus_p {\cal F}^p {\rm D}_{\rm crys}(V) \otimes_{K_0} \B_{\rm crys} \cap W_{2p}V\otimes_{\Q_p} \B_{\rm crys}.\\
\end{split}
\ee
On the other hand there is an isomorphism
\be \la{95}
\begin{split}
&  {\cal F}^p {\rm D}_{\rm crys}(V) \otimes_{K_0} \B_{\rm crys} \cap W_{2p}V\otimes_{\Q_p} \B_{\rm crys} \overset{\sim}{\rightarrow} \gr_{2p}^W V \otimes_{\Q_p} \B_{\rm crys}.\\
\end{split}
\ee Composing  isomorphism (\ref{94}) with the sum over $p$ of isomorphisms (\ref{95}) we get a canonical isomorphism:
\[\beta_{T}:\ V \otimes_{\Q_p} \B_{\rm crys}\overset{\sim}\rightarrow \bigoplus_p  \gr_{2p}^W V \otimes_{\Q_p} \B_{\rm crys}.\]
On the other hand a choice of a splitting $s$ of the weight filtration on $V$ provides an isomorphism
\[
\beta_{W, s}:\ \bigoplus_p \gr_{2p}^W V \otimes_{\Q_p} \B_{\rm crys} \overset{\sim}\rightarrow V \otimes_{\Q_p} \B_{\rm crys}.\]
The  composition of these isomorphisms is an endomorphism 
\[
\beta_{T} \circ \beta_{W, s} \in \End\Bigl(\bigoplus_p \gr_{2p}^W V \otimes_{\Q_p} \B_{\rm crys}\Bigr).
\]
Just like in the Hodge case, we call it   the p-adic crystalline period operator.  Its crucial for us property is that its matrix coefficients belong to the algebra ${\rm R}_{\rm crys}$. 
\[
\beta_{T} \circ \beta_{W, s} \in \End\Bigl(\bigoplus_p \gr_{2p}^W V \otimes_{\Q_p} {\rm R}_{\rm crys}\Bigr).
\]

\paragraph{The crystalline period morphism.}  Let $\widetilde {\cal V}$ be the set of split framed mixed  Tate crystalline    $\G_K$-representations. We define first a map of sets
\[
\Phi_{\rm crys}:  \widetilde {\cal V}    \longrightarrow {\rm T}({\rm B}_{\rm crys}^{\varphi =p}).
\]

For any split mixed  Tate crystalline representation $(V,s)$ and an $(i,j)$-framing 
$$
v \in \Hom(\mathbb{Q}_p(i), \gr^W_{-2i}V), ~f\in \Hom(\gr^W_{-2j}V, \mathbb{Q}_p(j))
$$ we define the period
\be
\label{period}
 \langle f\mid V, s \mid v\rangle := f\circ \beta_{T}\circ \beta_{W, s} \circ v \in \Hom_{{\B_{\rm crys}}}(\mathbb{Q}_p(i)\otimes {\B_{\rm crys}}, \mathbb{Q}_p(i+k) \otimes {\B_{\rm crys}})\cong {\B^{\varphi = p^{k}}_{\rm crys}}.
\ee
The last isomorphism is normalized so that it sends the map $  t^i \mapsto  t^{i+k}$ to $t^{k}$.

For each $m$, choose a basis $ \{v^\alpha_{m}\}$ in $\Hom(\mathbb{Q}_p(m), \gr^W_{-2m}V)$. Let $\{f_\alpha^{m}\}$ be the dual basis.

\bd
  \la{DEF4.1a}

  We define a map $\Phi_{\rm crys}$ recursively:
 \begin{equation}
  \la{REC1}
  \begin{split}
  &\Phi_{\rm crys}(f^n, V,v_0; s):=\\
  &-\sum_{0\leq m\leq n-1}  \sum_{\alpha}
  \langle f^n \mid V,s \mid v^\alpha_m \rangle t^{-(n-m-1)} \otimes t^{\otimes(n-m-1)} \otimes \Phi_{\rm crys}(f_\alpha^m, V, v_0;s).\\
  \end{split}
\end{equation}
\ed
 Definition \ref{DEF4.1a}  evidently does not depend on the choice of the generator $t$ of the Tate line.   
 It is   independent of the choice of bases. Elaborating (\ref{REC1}), we get:
\be
  \begin{split}
  &    \Phi_{\rm crys}(f^n, V,v_0; s):=\\
&    \sum_{\substack{1\leq k\leq n \\ 0=i_0<i_1<\cdots <i_k=n}}
     (-1)^k \bigotimes_{l=1}^{k} \left(\langle f^{i_l}\mid V,s \mid v_{i_{l-1}}\rangle t^{-(i_l-i_{l-1}-1)} \otimes t^{\otimes (i_l-i_{l-1}-1)}\right).\\
 \end{split}     
\ee 
 We can rewrite the definition of $\Phi_{\rm crys}$  as  
\begin{equation}
  \la{REC2}
  \sum_{0\leq m\leq n} \langle f^n \mid V,s \mid v_m\rangle  t^{-(n-m-1)} \otimes t^{\otimes(n-m-1)} \otimes \Phi_{\rm crys}(f^m, V, v_0;s)=0.
\end{equation}

 There is a natural projection 
$$
{\rm pr}: {\rm T}({\rm B}_{\rm crys}^{\varphi =p}) \lra {\cal A}^\varphi_\bullet({\rm R}_{\rm crys}).
$$
Let us set 
$$
{\mathcal{P}}_{\rm \G_K, crys}:={\rm pr}\circ \Phi_{\rm crys}.
$$ 
Just like in the Hodge-Tate case, one shows that the map ${\mathcal{P}}_{\rm \G_K, crys}$ does not depend on the splitting of the weight filtration, and that 
it sends equivalent framed mixed Tate Galois representations to the same element. So we get a well defined map of sets
\be \la{BCRYS}
\mathcal{P}_{\rm \G_K, crys}:\ {\cal H}_\bullet({\cal M}_{\rm \G_K, crys}) \rightarrow {\cal A}^\varphi_\bullet({\rm R}_{\rm crys}).
\ee

\bt \la{MAPPIxx}
The map $\mathcal{P}_{\rm \G_K, crys}$ is a homomorphism of graded Hopf algebras.
\et

\begin{proof} Follows literally the proof of Theorem \ref{MAPPI}, replacing everywhere $2\pi i$ by $t$.\end{proof}

\paragraph{The Hopf algebra ${\cal A}^\varphi_\bullet({\rm B}_{\rm crys}) $ and mixed Tate BKF-modules.}  The homomorphism (\ref{BCRYS}) is far from being an isomorphism. 
So a natural question arises, how to give an alternative definition for the  category of graded comodules over the Hopf algebra ${\cal A}^\varphi_\bullet({\rm B}_{\rm crys}) $?

Recall that L. Fargues \cite{F} introduced  a category,   called   now 
the category of Brueil-Kisin-Fargues modules. A variant, called the category of  {\it rigidified Brueil-Kisin-Fargues}   modules and denoted below by BKF, was 
considered in \cite{BMS}  and 
studied in \cite{A17}. The category BKF contains mutually non-isomorphic  invertible
 objects $A(n)$, $n \in \Z$, which are the analogs of the Tate motives, and one has \cite{A17}
\be \la{BKF}
\begin{split}
&{\rm Hom}_{{\cal M}_{\rm BKF }}({\rm A}(n), {\rm A}(n)) = \Q_p,\\
&{\rm Ext}_{\rm BKF}^1({\rm A}(0), {\rm A}(n)) = 
{\rm B}_{\rm dR}/t^n{\rm B}^+_{\rm dR}, ~~~~\forall n \in \Z. \\
\end{split}
\ee

Observe  that for $n > 0$  the ${\rm Ext}^1$-groups (\ref{CEG}) in the category ${\cal M}({\rm R}_{\rm crys})$ coincide  with the ${\rm Ext}^1$-groups (\ref{BKF})  in the category of 
rigidified Brueil-Kisin-Fargues   modules.

On the other hand  for $n \leq 0$ the ${\rm Ext}^1$-groups in the category ${\cal M}({\rm R}_{\rm crys})$  are zero, while the ones (\ref{BKF})  are not. 
To treat this problem we   consider  the category of mixed   BKF-modules, given by BKF-modules equipped with a weight filtration. 
Below we spell the definition of its subcategory of mixed Tate BKF-modules. 

\bd The category ${\cal M}_{\rm T- BKF}$ of mixed Tate ${\rm BKF}$-modules consists of ${\rm BKF}$-modules $M$ equipped with an increasing weight filtration $W_\bullet$ such that 
for each integer $n$ 
the associate graded ${\rm gr}^W_{-2n}M$ is a direct sum of copies of ${\rm A}(n)$ and ${\rm gr}^W_{-2n+1}M=0$.  The morphisms are the ones in the original category which are strongly compatible with the weight filtration. 
\ed

Then the objects $A(n)$ are mutually non-isomorphic, and one has 
\be \la{BKF1}
\begin{split}
&{\rm Hom}_{{\cal M}_{\rm T-BKF }}({\rm A}(m), {\rm A}(m)) = \Q_p, ~~~~\forall m \in \Z,\\
&{\rm Ext}_{{\cal M}_{\rm T-BKF }}^1({\rm A}(0), {\rm A}(n)) = {\rm B}_{\rm dR}/t^n{\rm B}^+_{\rm dR}, ~~~~  n \in \Z_{>0}.\\
&{\rm Ext}_{{\cal M}_{\rm T-BKF }}^1({\rm A}(0), {\rm A}(n)) = 0, ~~~~  n \in \Z_{\leq 0}.\\\end{split}
\ee
 
\bl \la{LEMM5.3}
If ${\rm Ext}_{\rm BKF}^{>1}({\rm A}(0), {\rm A}(n)) =0$ for all $n\in \Z$, then the category ${\cal M}({\rm R}_{\rm crys})$
is equivalent to the category ${{\cal M}_{\rm T-BKF}}$ of  mixed Tate ${\rm BKF}$-modules. 
\el 

\begin{proof} Follows immediately from general properties of mixed Tate categories, combined with the fact that their Ext-groups coincide. \end{proof}

Equivalently, under the assumption of Lemma \ref{LEMM5.3} the fundamental Hopf algebra ${\cal H}_\bullet({{\cal M}_{\rm T-BKF}})$ is isomorphic to 
our graded Hopf algebra ${\cal A}^\varphi_\bullet({\rm R}_{\rm crys})$.   

Conjecture \ref{CONN5.2} below just means  that there exists a \underline{canonical} period isomorphism 
$$
{\cal P}_{\rm crys}: {\cal H}_\bullet({{\cal M}_{\rm T-BKF}}) \stackrel{\sim}{\lra} {\cal A}^\varphi_\bullet({\rm R}_{\rm crys})
 $$
 which fits into a commutative diagram 
 
   \centerline{
    \xymatrix{
     {\cal H}_\bullet({\cal M}_T(F, v))  \ar[rr]^{\mbox{\rm p-adic realisation}} \ar[rd]_{{\cal P}_{\rm T, crys}}& &{\cal H}_\bullet({{\cal M}_{\rm T-BKF}})  \ar[ld]_{\sim}^{{\cal P}_{\rm crys}}\\
       & {\cal A}^\varphi_\bullet({\rm R}_{\rm crys})&  }  }

\bcon \la{CONN5.2}
There is a \underline{canonical} equivalence between the  category  category   ${{\cal M}_{\rm T-BKF}}$ of  mixed Tate ${\rm BKF}$-modules and the category 
${\cal M}({\rm B}_{\rm crys})$ of graded ${\cal A}^\varphi_\bullet({\rm R}_{\rm crys})$-comodules:
$$
{\rm P}:  {{\cal M}_{\rm T-BKF}} \stackrel{\sim}{\lra} {\cal M}({\rm R}_{\rm crys}).
$$
\econ

We will construct the map  ${\cal P}_{\rm crys}$ and prove Conjecture \ref{CONN5.2} in \cite{GZ}.

\paragraph{The   crystalline period map for  mixed Tate motives over a number field. } Let $F$ be a  number field. Let $v$ be a   
  non-archimedian  place  of $F$ over $p$.   Denote by  ${\cal M}_T(F, v)$  the subcategory of the category of mixed Tate motives 
  ${\cal M}_T(F)$ over $F$  unramified at 
  $v$  \cite{DG}.  It is a mixed Tate category. So there is   
    the fundamental Hopf algebra ${\cal H}_\bullet({\cal M}_T(F, v))$ of this category. 
    
    Denote by $F_v$ the completion of $F$ at $v$.  
Let us assign to a mixed Tate motive $M$ over $F$ the corresponding p-adic  Galois representation of $\G_{F_v}$, which is 
   evidently mixed Tate, and crystalline if   $M$ is unramified at $v$. So we arrive at the   p-adic realization  functor  
\be \la{CANRF}
   {\cal M}_T(F, v) \lra {\cal M}_{\rm \G_{F_v}, crys}.
\ee
      
 The functor (\ref{CANRF}) is a  functor between mixed Tate categories.    The fundamental Hopf algebra of a mixed Tate category can be constructed  
 via the set of equivalence classes of  framed objects in the   category. 
      Thus the p-adic realization  functor gives rise to a map of graded Hopf algebras
  $$
  {\cal R}_{F_v}: {\cal H}_\bullet({\cal M}_T(F, v)) \lra {\cal H}_\bullet({\cal M}_{\rm \G_{F_v}, crys}).
  $$
Composing it with the map provided by Theorem \ref{MAPPIxx} we get  a map of graded Hopf algebras
  $$
  {\cal P}_{\rm T, crys}: =  \mathcal{P}_{\rm \G_K, crys}\circ {\cal R}_{F_v}: {\cal H}_\bullet({\cal M}_T(F, v))\lra  {\cal A}^\varphi_\bullet({\rm R}_{\rm crys}).
     $$

\paragraph{The p-adic regulator map.} 
\bc Let $F$ be a number field. Then there is a p-adic regulator map 
\be \la{BKF1as}
\begin{split}
&K_{2n-1}(F)\lra  {\rm B}_{\rm dR}/t^n{\rm B}^+_{\rm dR}, ~~~~  n \in \Z_{>0}.\\
\end{split}
\ee
\ec 
 
\begin{proof} For a number field $F$ one has the basic isomorphism 
$$
K_{2n-1}(F) \otimes \Q = {\rm Ext}_{{\cal M}_{\rm T}(F)}^1({\Q}(0), {\Q}(n)).
$$
So combining this isomorphism with the   isomorphism (\ref{CEG})   we arrive at the map (\ref{BKF1as}) for $n>1$.

 For $n=1$ we need to use the semi-stable story,  discussed in Section \ref{sec5.4}. 
\end{proof}

\subsection{The period morphism ${\cal P}_{\rm st}$ for semi-stable mixed Tate Galois modules} \la{sec5.4}

  \paragraph{Fontaine's semi-stable period algebra $\B_{\rm st}$.} The  algebra $\B_{\rm st}$ is a $\G_K$-module which contains the algebra $\B_{\rm crys}$. There is a $\G_K$-equivariant
  isomorphism of algebras  with $u="\log [p]"$:
  $$
 \B_{\rm st} = \B_{\rm crys}[u], ~~~~g\cdot u = u + a(g)  t, ~~\forall g \in \G_K.
 $$
The ring $\B_{\rm st}$ is equipped with the Frobenious map $\varphi$ and the monodromy operator 
$N$, 
such that  
$$
\varphi (u) = p u, ~~~~N = -\frac{d}{du}, ~~~~N\varphi = p \varphi N.
$$ 
Evidently, there is an exact sequence of vector spaces
$$
0 \lra \B_{\rm crys} \lra \B_{\rm st}  \stackrel{N}{\lra} \B_{\rm st} \lra 0.
$$

Consider the $\Q_p$-vector space 
$$
{\rm D}_{\rm st}(V) := \Bigl(V \otimes_{\Q_p}  {\rm B}_{\rm st}\Bigr)^{\G_{K}}.
$$
A p-adic $\G_K$-representation $V$ is called {\it semi-stable} if  ${\rm dim}{\rm D}_{\rm st}(V) = {\rm dim}(V)$.

\paragraph{The algebra ${\rm R}_{\rm st}$.} The algebra $\B_{\rm st}$ has a subalgebra 
$$
{\rm R}_{\rm st} : = \bigoplus_{m \in \Z} {\rm B}^{\varphi = p^m}_{\rm st}.
 $$  
The monodromy operator $N$ restricts to the subalgebra ${\rm R}_{\rm st}$. Indeed,
$$
N: \B^{\varphi = p^m}_{\rm st} \lra \B^{\varphi = p^{m-1}}_{\rm st}.
$$
Evidently, there is an exact sequence of vector spaces
$$
0 \lra {\rm R}_{\rm crys} \lra {\rm R}_{\rm st}  \stackrel{N}{\lra} {\rm R}_{\rm st} \lra 0.
$$

\paragraph{A dg-algerbra ${\rm R}^\bullet_{\rm st}$.}
We interpret the monodromy operator $N$ as a differential in the complex ${\Bbb B}^\bullet_{\rm st}:= \B_{\rm st}  \stackrel{N}{\lra} \B_{\rm st}$, placed  
   in  degrees $[0,1]$. 
The commutative algebra structure on  $\B_{\rm st}$ extends to a commutative
dg-algebra structure on the complex ${\Bbb B}^\bullet_{\rm st}$  as follows. 
For any $b_i, b_i' \in {\Bbb B}^i_{\rm st}$ set
\be
\mu(b_0\otimes b'_0) := b_0b'_0, ~~~~\mu(b_0\otimes b_1) = \mu(b_1\otimes b_0):= b_0b_1, ~~~~ \mu(b_1\otimes b'_1) =0. 
\ee
The differential  satisfies the Leibniz rule. 

There is a dg-subalgebra ${\rm R}^\bullet_{\rm st}$ of the dg-algebra ${\Bbb B}^\bullet_{\rm st}$:
 \be
\begin{split}
& {\rm R}^0_{\rm st} \stackrel{d}{\lra}  {\rm R}^1_{\rm st}\\
&{\rm R}^0_{\rm st}= {\rm R}^1_{\rm st}:= {\rm R}_{\rm st}, ~~~~d:= N.\\
  \end{split}
\ee 
There is a $\varphi$-grading on the complex ${\rm R}^\bullet_{\rm st}$:
$$
{\rm R}^0_{\rm st} = \bigoplus_{m \in \Z} {\rm R}^{0, \varphi = p^m}_{\rm st}, ~~~~ {\rm R}^1_{\rm st} = \bigoplus_{m \in \Z} {\rm R}^{0, \varphi = p^{m-1}}_{\rm st},~~~~{\rm deg}_\varphi(d) = -1.
$$

Let us axiomatize the properties of the   dg-algebra ${\rm R}^\bullet_{\rm st}$.

 \bd A  Tate $\varphi$-dg algebra $R= \oplus_{m,n\in \Z}R^n_m$   is  a dg--algebra  over a field $k$  equipped with an additional 
 $\varphi$-grading, 
   a differential $d: R_m^{n} \lra  R_{m-1}^{n+1}$ of the cohomological degree $+1$ and  the $\varphi$-degree  $  -1$, and  
 an invertible Tate line  $k(1)  \subset R^0_1$, as in (\ref{invtl}).
 \ed

   \paragraph{The Hopf dg-algebra associated to a Tate $\varphi$-dg algebra.}
   
    A  Tate $\varphi$-dg algebra $R$ gives rise to a graded $k$-vector space  ${\cal A}^\varphi_\bullet(R)$:
 \be \la{phial}
\begin{split}
  &{\cal A}^\varphi_\bullet(R):=   \bigoplus_{n=0}^\infty {\cal A}^\varphi_n(R), ~~~~{\cal A}^\varphi_0(R) = k.\\
  &{\cal A}^\varphi_d(R):=        \overline {R_1 } \otimes_{k} \underbrace{{R_1 } \otimes_{k}
\ldots \otimes_{k} {R_1 }}_{\text{ $d-1$ factors }}, ~~d>0.\\
\end{split}
\ee

 The differential $d$ in $R$ provides it with a differential 
 \be \la{DIFFDDab}
{\cal D}(\overline x_1 \otimes x_2 \otimes \ldots \otimes x_n):= \overline {d\overline x_1 \cdot  x_2} \otimes (x_3 \otimes \ldots \otimes x_n) + (-1)^{|x_1|} \overline x_1 \otimes {\cal D}(x_2 \otimes \ldots \otimes x_n).
\ee 
  It is well defined since $x,y  \in R_1$ and $d$ has $\varphi$-degree $-1$, so $d x \cdot  y  \in R_1$. So the same arguments as in the proofs of 
  Theorems \ref{TH2.2} and \ref{Th2.2v} deliver the following basic result.

 \bt \la{Th2.2vv} Let $R$ be a   Tate $\varphi$-dg algebra. Then the map ${\cal D}$ is a differential, providing 
  ${\cal A}^\varphi_\bullet(R)$ with a   Hopf dg-algebra structure.  
\et

  \paragraph{Hopf dg-algebra ${\cal A}^\varphi_\bullet({\rm R}^\bullet_{\rm st})$.} 
By Theorem \ref{Th2.2vv} the Tate $\varphi$-dg algebra ${\rm R}^\bullet_{\rm st}$ gives rise to a   Hopf  dg algebra ${\cal A}^\varphi_\bullet({\rm R}^\bullet_{\rm st})$:
 \be
\begin{split}
&{\cal A}^\varphi_\bullet({\rm R}^\bullet_{\rm st}) = \Q_p \oplus \bigoplus_{n=1}^\infty{\cal A}^\varphi_n({\rm R}^\bullet_{\rm st}),\\
&{\cal A}^\varphi_n({\rm R}^\bullet_{\rm st}):=    
  {\rm B}_{\rm dR}/t{\rm B}^+_{\rm dR} \otimes_{\Q_p} \underbrace{{\rm B}_{\rm st}^{\varphi =p}
 \otimes _{\Q_p}   \ldots \otimes _{\Q_p}  {\rm B}_{\rm st}^{\varphi =p}}_{\text{ $n-1$ factors }}.\\
\end{split}
\ee 
   
Taking $H^0_{\cal D}$, we get a graded Hopf algebra $H^0_{\cal D}({\cal A}^\varphi_\bullet({\rm R}^\bullet_{\rm st}))$

   \paragraph{The p-adic  period operator.}

Let $V$ be a mixed  Tate semi-stable $\G_K$-representation over $\Q_p$.   Just as in the crystalline case, there is a decreasing slope filtration on 
${\cal F}^\bullet$ on $ {\rm D}_{\rm crys}(V)$, and an increasing weight filtration $W_\bullet$ on $V$. 
Thanks to (\ref{91}),   the slope filtration 
${\cal F}^p {\rm D}_{\rm st}(V) \otimes_{\Q_p} \B_{\rm st}$ is   opposite to the weight filtration $W_\bullet V \otimes_{K_0} \B_{\rm st}$. 
Therefore these two filtrations split $V \otimes_{\Q_p} \B_{\rm st}$. So there are two isomorphisms: 
\be
\begin{split}
&V \otimes_{\Q_p} \B_{\rm st} =\bigoplus_p {\cal F}^p {\rm D}_{\rm st}(V) \otimes_{K_0} \B_{\rm st} \cap W_{2p}V\otimes_{\Q_p} \B_{\rm st},\\
&  {\cal F}^p {\rm D}_{\rm st}(V) \otimes_{K_0} \B_{\rm st} \cap W_{2p}V\otimes_{\Q_p} \B_{\rm st} \overset{\sim}{\rightarrow} \gr_{2p}^W V \otimes_{\Q_p} \B_{\rm st}.\\
\end{split}
\ee 
Combining them,  we get a canonical isomorphism:
\[\beta_{T}:\ V \otimes_{\Q_p} \B_{\rm st}\overset{\sim}\rightarrow \bigoplus_p  \gr_{2p}^W V \otimes_{\Q_p} \B_{\rm st}.\]
A choice of a splitting $s$ of the weight filtration on $V$ provides an isomorphism
\[
\beta_{W, s}:\ \bigoplus_p \gr_{2p}^W V \otimes_{\Q_p} \B_{\rm st} \overset{\sim}\rightarrow V \otimes_{\Q_p} \B_{\rm crys}.\]
Let us consider their composition,  which we call    the p-adic semi-stable period operator. Since we are in the mixed Tate case, its matrix coefficients are in the subalgebra 
${\rm R}_{\rm st}$: 
\[
\beta_{T} \circ \beta_{W, s} \in \End\Bigl(\bigoplus_p \gr_{2p}^W V \otimes_{\Q_p} {\rm R}_{\rm st}\Bigr).
\]

Repeating verbatim the construction of the map $\Phi$ in the semi-stable set-up,  we get a homomorphism of graded Hopf algebras:
$$
\mathcal{P}_{\rm \G_K, st}:\ {\cal H}_\bullet({\cal M}_{\rm \G_K, st}) \longrightarrow {\cal A}^\varphi_\bullet({\rm R}_{\rm st}).
$$
Unlike the crystalline case,  in the semi-stable case there is an extra condition the image of the  morphism $\mathcal{P}_{\rm \G_K, st}$ satisfies. To state
 it, observe that 
there is an inclusion of Hopf dg-algebras: 
$$
H^0_{\cal D}({\cal A}^\varphi_\bullet({\rm R}^\bullet_{\rm st}))  \subset {\cal A}^\varphi_\bullet({\rm R}^\bullet_{\rm st}).
$$
Indeed, by  definition the cohomological degree $0$ part of ${\cal A}^\varphi_\bullet({\rm R}^\bullet_{\rm st})$ 
  coincides with  ${\cal A}^\varphi_\bullet({\rm R}_{\rm st})$.  So
 $$
 H^0_{\cal D}({\cal A}^\varphi_\bullet({\rm R}^\bullet_{\rm st})) = {\rm Ker}({\cal D}) \cap  {\cal A}^\varphi_\bullet({\rm R}_{\rm st}).
 $$

\bt \la{MAPPIxx}
The image of the map $\mathcal{P}_{\rm \G_K, st}$ lies in the  graded Hopf subalgebra $H^0_{\cal D}({\cal A}^\varphi_\bullet({\rm R}^\bullet_{\rm st}))$:
$$
\mathcal{P}_{\rm \G_K, st}:\ {\cal H}_\bullet({\cal M}_{\rm \G_K, st}) \longrightarrow H^0_{\cal D}({\cal A}^\varphi_\bullet({\rm R}^\bullet_{\rm st})).
$$
\et

\begin{proof} The proof follows the lines of the proof of Theorem \ref{MAPPP2}. Details will appear in \cite{GZ}.
\end{proof}

The dg-category of dg-comodules 
over the dg-Hopf algebra ${\cal A}^\varphi_\bullet({\rm R}^\bullet_{\rm st})$ is a dg-model for the derived category of the abelian category of p-adic mixed 
Tate semi-stable $\G_K$-representations: 

\bt There is a canonical equivalence between the derived category of the abelian category of mixed Tate semi-stable $\G_K$-representations over $\Q_p$ 
and the $H^0$ of the dg-category of dg-comodules over $\Q_p$ 
over the dg-Hopf algebra ${\cal A}^\varphi_\bullet({\rm R}^\bullet_{\rm st})$.
\et

  The details of the proof of this Theorem will appear in \cite{GZ}.

\section{Appendix A: Mixed Tate categories} \la{SSEc5}

 Our exposition  borrows from  Section 1 of the unpublished work \cite{BD92}.

Let ${\cal T}$ be a Tannakian $F$-category, where $F$ is a field, and let $F(1)$  be a fixed rank one object of ${\cal T}$. A pair $({\cal T}, F(1))$ is called a {\it mixed Tate category}
if the objects
$F(n):=F(1)^{\otimes n}$, $n \in \Z$, are mutually non-isomorphic, any irreducible object in ${\cal T}$ is isomorphic to some $F(i)$, and one has
$$
{\rm Ext}^1_{\cal T}({F}(0), F(n)) =0, ~~~~  \forall n \leq 0.
$$

Let $({\cal T}_1, {\cal T}_2)$ be a pair of mixed Tate categories. A {\it pure functor} $\rho:{\cal T}_1\lra  {\cal T}_2$ is a tensor functor equipped with an isomorphism
$\rho(F(1)_{{\cal T}_1}) =  F(1)_{{\cal T}_2}$.

  Lemma \ref{BD0} establishes a canonical weight filtration on the objects of a mixed Tate category.

\bl \la{BD0} Any object $M$ in a mixed Tate category ${\cal T}$ has a unique increasing finite filtration $W_\bullet$,   indexed by even integers, such that ${\rm gr}^W_{2k}M$ is a direct sum of copies of $F(-k)$.
Any morphism is strictly compatible with the weight filtration $W_\bullet$.
\el

Given a mixed Tate category $({\cal T}, F(1))$,   there is a fiber functor to the category of finite dimensional graded $F$-vector spaces, called   {\it the canonical fiber functor}:
$$
\omega_{\cal T}:
{\cal T} \lra {\rm Vect}_F^\bullet, \qquad X\lra \oplus_n{\rm Hom}({F}(-n), {\rm gr}^W_{2n}X).
$$
Indeed,  it is an exact tensor functor by Lemma \ref{BD0}. It is a pure functor.

Denote by ${\rm L}_{\cal T}$ the graded pro-Lie algebra of $\otimes$-derivations of the functor $\omega_{\cal T}$:
$$
{\rm L}_{{\cal T}}:= {\rm Der}^\otimes(\omega_{\cal T}).
$$
An element $\alpha \in {\rm L}_{{\cal T}, i}$ is a degree $i$ endomorphism of $\omega_{\cal T}$, that is a collection of morphisms of
 functors $\alpha_i: \omega_j \to \omega_{j+i}$ such that, using   isomorphisms $\omega(M\otimes N) =  \omega(M ) \otimes \omega( N)$ and $\omega(M(1) )_j = \omega(M)_{j+1}$,
  we have:
  $$
\alpha_{M\otimes N} = \alpha_{M } \otimes {\rm Id}_{\omega(N)} + {\rm Id}_{\omega(M)}  \otimes \alpha_{ N}, ~~~~\alpha_{j, M(1)} = \alpha_{j+1, M} .
$$
Since the morphisms are strictly compatible with the weight filtration, we have ${\rm L}_{{\cal T}, i}=0$ if $i \geq 0$.

The  graded Lie algebra ${\rm L}_{{\cal T}}$ is called the {\it fundamental Lie algebra} of the mixed Tate category ${\cal T}$.
It is isomorphic, although non-canonically, to a graded Lie algebra generated by the graded vector space
$$
 \mbox{ ${\rm Ext}^1_{\cal T}(F(0), F(n))^*$ in the degree $-n$, where $n = 1,2, ...$},
$$
subject to the quadratic relations between the generators given in the degree $-n$ by the vector space ${\rm Ext}^2_{\cal T}(F(0), F(n))^*$.
Precisely, the dual to the product of the Ext's provides us a map
$$
{\rm Ext}^2_{\cal T}(F(0), F(n))^* \lra \bigoplus_{p+q=n} {\rm Ext}^1_{\cal T}(F(0), F(p))^* \wedge{\rm Ext}^1_{\cal T}(F(0), F(q))^*.
$$
The space of the degree $-n$ relations is the image of this map.

Given any negatively graded Lie algebra $L$ over a field $F$, the category ${\cal R}(L)$ of graded finite dimensional representations of $L$ is a mixed Tate category, with
the object $F(1)$ given by the trivial one dimensional representation   in the degree $-1$.
The functor
$\omega$ is an equivalence of  mixed Tate categories
\be \la{eqcc}
{\cal T} \stackrel{\sim}{\lra} {\cal R}({\rm L}_{{\cal T}}).
\ee

\bd The fundamental Hopf algebra ${\cal H}_\bullet({\cal T})$ of a mixed Tate category ${\cal T}$ is the graded linear dual to the 
universal enveloping algebra of the fundamental Lie algebra ${\rm L}_{{\cal T}}$.  It is a commutative graded Hopf algebra. 
\ed

The mixed Tate category ${\cal T}$ is canonically equivalent to the category of finite dimensional graded comodules over the  fundamental Hopf algebra ${\cal H}_\bullet({\cal T})$.

\section{Appendix B: Quasi-shuffle product  of multiple  polylogarithms}
The  formula for the product $m'$ looks identical to the
  product formula for the multiple polylogarithm series, since both are defined using the quasi-shuffle product. Let us elaborate on this.
  Recall   \cite{G94} that the depth $m$ multiple polylogarithm  power series\footnote{In fact the definition
we use differs slightly from the one in {\it loc. cit.} by using the non-strict inequalities $... k_i \leq k_{i+1} ... $ instead of the strict ones.
The advantage of the latter is  that the function defined by the resulting power is equal to a single iterated integral \cite[Theorem 16]{G94}.}
are determined by a collection of positive integers $p_1, ..., p_m$:
$$
  {\rm Li}^\prime_{p_1, ..., p_m}(x_1, ..., x_m):= \sum_{0<  k_1 \leq \ldots \leq k_m}\frac{x_1^{k_1} \ldots x_m^{k_m}}{k^{p_1}_1 \ldots k^{p_m}_m}.
$$
The integer $w = s_1 + \ldots + s_m$ is called the weight.
The product
$$
{\rm Li}^\prime_{p_1, ..., p_m}(x_1, ..., x_m) \cdot {\rm Li}^\prime_{p'_1, ..., p'_n}(x'_1, ..., x'_n)
$$
of two multiple polylogarithm power series of weights $w$ and $w'$ is  a sum  of multiple polylogarithm series of the weight $w+w'$. For example we have
 --- compare with the first
line in (\ref{QSPFa}) ---
\be
\begin{split}
&{\rm Li}^\prime_{1}(x) \cdot {\rm Li}_{1}^\prime(y) = {\rm Li}^\prime_{1, 1}(x,y) + {\rm Li}^\prime_{1, 1}(y,x) - {\rm Li}_{2}^\prime(xy).\\
\end{split}
\ee
Indeed, just as in (\ref{QS1}), the   product sum is  decomposed into a disjoint union of three subsums:
$$
\{ k_1\}\cup \{k_2\} = \{ k_1 \leq k_2\} \cup \{ k_2 \leq k_1\}  - \{k_1=k_2\}.
$$
Similarly, we have ---  compare with the second
formula in (\ref{QSPFa}) ---
$$
 {\rm Li}^\prime_{1}(x) \cdot {\rm Li}^\prime_{1,1}(y,z) = {\rm Li}^\prime_{1, 1,1}(x,y,z) +  {\rm Li}^\prime_{1, 1,1}(y,x,z) +  {\rm Li}^\prime_{1, 1,1}(y,z,x) - {\rm Li}^\prime_{2,1}(xy, z) -{\rm Li}^\prime_{1, 2}(y, xz).
$$

We assign, entirely   \underline{formally}, to  weight $w$ multiple polylogarithms  elements of  $\C^{\otimes w}$:
\be \la{FORC}
\bs
& \ {\rm Li}_{p_1, ..., p_m}(x_1, ..., x_m) \ \longmapsto \  \\
&\frac{x_1}{(2\pi i)^{p_1-1} } \otimes (2\pi i)^{\otimes p_1-1} \otimes \frac{x_2}{(2\pi i)^{s_2-1}} \otimes (2\pi i)^{\otimes p_2-1} \otimes \ldots \otimes \frac{x_m }{(2\pi i)^{p_m-1}}\otimes (2\pi i)^{\otimes s_m-1}.\\
\end{split}
\ee

\bl \la{T1.6}
The formal map (\ref{FORC}) transforms the product formula for the multiple polylogarithm power series into the product $m'$ in ${\rm T}(\C)$.
\el

This follows easily from Lemma \ref{L4.4}.
Let us give one more example:
$$
{\rm Li}_{p}(x) \cdot {\rm Li}_{q}(y) = {\rm Li}_{p,q}(x,y) + {\rm Li}_{q,p}(y,x) - {\rm Li}_{p+q}(xy).
$$
In complete analogy with this, one   verifies  using   inductive definition (\ref{EQP}), see Lemma \ref{L4.4},  that
\be
\begin{split}
&\Bigl(\frac{x}{(2\pi i)^{p-1}} \otimes (2\pi i)^{p-1}\Bigr)\ast_{m'} \Bigl(\frac{y}{(2\pi i)^{q-1}} \otimes (2\pi i)^{q-1}\Bigr) = \\
& \frac{x}{(2\pi i)^{p-1}}   \otimes (2\pi i)^{p-1}\otimes y \otimes (2\pi i)^{q-1} + \frac{y}{(2\pi i)^{q-1}}  \otimes (2\pi i)^{q-1}  \otimes \frac{x}{(2\pi i)^{p-1}}  \otimes (2\pi i)^{p-1}\\
& - \frac{xy}{(2\pi i)^{p+q-1}}  \otimes (2\pi i)^{p+q-1}.\\
 \end{split}
 \ee

\paragraph{Warning.} Multiple polylogarithms are periods of mixed Hodge-Tate structures, and thus give rise to  elements
of the Hopf algebra ${\cal H}_\bullet$. However the map (\ref{FORC}) \underline{does not} respect the coproducts.

\end{document}